\documentclass[a4paper]{article}

\usepackage{amsmath,amssymb,amscd, amsthm}

\newtheorem*{thm}{Theorem}
\newtheorem*{prop}{Proposition}
\newtheorem*{lem}{Lemma}

\newtheorem*{dfn}{Definition}
\newtheorem*{cor}{Corollary}
\newtheorem*{exa}{Example}
\newtheorem*{rem}{Remark}

\newcommand{\Z}{\mathbb{Z}}
\newcommand{\C}{\mathbb{C}}
\newcommand{\R}{\mathcal{R}}

\newcommand{\M}{\mathfrak{M}}
\newcommand{\qZ}{\mathfrak{Z}}

\newcommand{\tL}{\mathfrak{L}}
\newcommand{\B}{\mathcal{B}}

\newcommand{\dimv}{\mathbf{v}}
\newcommand{\dimw}{\mathbf{w}}

\newcommand{\maya}{\mathbf{k}}

\newcommand{\hilb}{\left(\C^2\right)^{[n]}}

\newcommand{\node}[1]{{\hspace{1pt}\text{\setlength{\fboxsep}{.5mm}\fbox{$#1$}}\hspace{1pt}}}

\newcommand{\bk}{\mathbf{k}}

\newcommand{\bm}{\mathbf{m}}
\newcommand{\bn}{\mathbf{n}}
\newcommand{\bv}{\mathbf{v}}
\newcommand{\bj}{\mathbf{j}}
\newcommand{\bal}{\mathbf{d}}
\newcommand{\alp}{d}

\newcommand{\ha}{\mathbf{H}_N}
\newcommand{\aha}{\dot{\mathbf{H}}_N}
\newcommand{\ahay}{\dot{\mathbf{H}}^{(1)}_N}
\newcommand{\ahax}{\dot{\mathbf{H}}^{(2)}_N}
\newcommand{\daha}{\ddot{\mathbf{H}}_N}
\newcommand{\poly}{\mathcal{R}[z_1^{\pm 1},\ldots,z_N^{\pm 1}]}

\newcommand{\tor}{U_{\R}'(\mathfrak{sl}_{\,l,tor})}

\newcommand{\uqd}{{U^{(1)}_\R}'(\hat{\mathfrak{sl}}_l)}
\newcommand{\uqc}{{U^{(2)}_\R}'(\hat{\mathfrak{sl}}_l)}
\newcommand{\uq}{U_q'(\hat{\mathfrak{sl}}_l)}

\newcommand{\vt}{\overset{v}{T}}
\newcommand{\pt}{\overset{p}{T}}

\newcommand{\ubm}{{\underline{\bm}}}
\newcommand{\ubj}{{\underline{\bj}}}
\newcommand{\um}[1]{{\underline{m}_{\,{#1}}}}
\newcommand{\uj}[1]{{\underline{j}_{\,{#1}}}}
\newcommand{\sj}{{\widehat{\sigma_\bj}}}
\newcommand{\sm}{{\widehat{\sigma_{\ubm}}}}

\title{K-theory of quiver varieties, q-Fock space\\and nonsymmetric Macdonald polynomials}
\author{Kentaro Nagao}

\begin{document}

\maketitle

\begin{abstract}
We have two constructions of the level-$(0,1)$ irreducible representation of the quantum toroidal algebra of type $A$.
One is due to Nakajima and Varagnolo-Vasserot. 
They constructed the representation on the direct sum of the equivariant K-groups of the quiver varieties of type $\hat{A}$.
The other is due to Saito-Takemura-Uglov and Varagnolo-Vasserot. 
They constructed the representation on the q-deformed Fock space introduced by Kashiwara-Miwa-Stern.

In this paper we give an explicit isomorphism between these two constructions.
For this purpose we construct simultaneous eigenvectors  on the q-Fock space using nonsymmetric Macdonald polynomials.
Then the isomorphism is given by corresponding these vectors to the torus fixed points on the quiver varieties.

\end{abstract}

\tableofcontents

\section{Introduction}

{\ }

\vspace{-12pt}

Geometry of \textbf{quiver varieties}, introduced by Nakajima, involves rich mathematical structures.

One of the most notable result is Nakajima's construction of representations of quantum loop algebras (\cite{quiver3}) : 
the direct sum of torus equivariant K-groups of the quiver varieties is endowed with a structure of a representation of 
the quantum affinization $U_\R(L\mathfrak{g})$ of the corresponding Kac-Moody algebra $\mathfrak{g}$.
The resulting representation is what we call an $l$-highest weight representation, 
that is to say, a "highest weight representation" with respect to the triangular decomposition of the  quantum loop algebra.
 
Let us concentrate our attention on quiver varieties of affine type. 
They appear in gauge theory as framed moduli spaces of instantons on ALE spaces, 
which originally motivated Nakajima to introduce quiver varieties.
They also have interesting connections with some areas in mathematics such as the theory of McKay correspondence and 
the representation theory of symplectic reflection algebras (see \cite{haiman-lec-note} and \cite{gordon-ohp} for example). 
In this point of view, more careful study about the actions of the quantum toroidal algebras, 
quantum affinizations of the affine Kac-Moody algebras, on the equivariant K-groups seems to be important.

\medskip

Schur-Weyl duality is an equivalence between certain categories of representations of $\mathfrak{gl}_l$ and of $\mathfrak{S}_n$. 
Varagnolo-Vasserot show that there exists an analogous duality between the quantum toroidal algebra of type $A$ and 
a certain double affinization of the Hecke algebra of type $A$, called the toroidal Hecke algebra (\cite{vv-schur}).  
The toroidal Hecke algebra has a remarkable representation called Dunkl-Cherednik representation.
Applying Schur-Weyl duality for Dunkl-Cherednik representation, 
Saito-Takemura-Uglov and Varagnolo-Vasserot construct the representation of the quantum toroidal algebra (\cite{saito-takemura-uglov}, \cite{vv-q-fock}). 
The underlying space is so called the \textbf{q-Fock space} (\cite{kashiwara-miwa-stern}).

\medskip

There are much fewer things known about representations of quantum toroidal algebras than of quantum affine algebras.
Now, at least, we have two constructions of the representation of the quantum toroidal algebra of type $A$.
In this paper we give an explicit isomorphism between these two constructions. 
We hope it will be helpful for further analyses of the representation, such as study of canonical bases of the representations.

\medskip

We can describe the representation on the equivariant K-groups in a combinatorial manner using the localization theorem (\cite{vv-k}). 
In particular, the torus fixed points correspond to simultaneous eigenvectors for the action of a certain subalgebra of the quantum toroidal algebra.
Our strategy is to construct simultaneous eigenvectors on the q-Fock space.
The isomorphism will be given by corresponding these vectors to the torus fixed points.
For the construction of simultaneous eigenvectors, \textbf{nonsymmetric Macdonald polynomials} plays a crucial role, 
where nonsymmetric Macdonald polynomials are simultaneous eigenvectors for Dunkl-Cherednik operators (\cite{cherednik-nsmp}, \cite{macdonald-bourbaki}, \cite{opdam-nsmp}),

Takemura-Uglov described the irreducible decomposition of the q-Fock space as the representation of a certain subalgebra of the quantum toroidal algebra, 
which is isomorphic to the quantum affine algebra (\cite{takemura-uglov-level0}).
They also showed that each irreducible components are isomorphic to tensor products of fundamental representations. 
For this purpose they introduced specific vectors of the q-Fock space using nonsymmetric Macdonald polynomials (see Remark \ref{eigenvec}).

In this paper we introduce new vectors.
We also use nonsymmetric Macdonald polynomials, but an additional operation is required (see \ref{renumbering}). 
They are simultaneous eigenvectors and the main subject of this paper.
These simultaneous eigenvectors 
allow us a combinatorial description of the representation on the q-Fock space 
and we can see this coincides with the combinatorial description of the representation on the equivariant K-groups.

\medskip

In \cite{vv-decomp} and \cite{schiffmann-hall-and-fock}, the action of the Hall algebra of the cyclic quiver on the q-Fock space is studied.
The Hall algebra of the cyclic quiver is realized using perverse sheaves on the space of representations of the quiver by Lusztig (\cite{lusztig-quiver}).
Nakajima's construction of quiver varieties and representations on their K-groups are, philosophically, parallel to Lusztig's construction.
We could expect this observation gives conceptual interpretation of the isomorphism constructed in this paper.     
In particular, this isomorphism would help us to study of canonical bases of the K-groups of quiver varieties (see \cite{lusztig-rem-on-quiver-var} and \cite{vv-canonical-basis-of-quiver-var}, for quiver varieties of finite type). 

\medskip

In \S \ref{k-theory} -- \S \ref{toroidal on q} we are mainly occupied with review of, and arrangement for our use of, the results of \cite{vv-k}, \cite{vv-schur}. \cite{saito-takemura-uglov} and \cite{vv-q-fock}. 
In \S \ref{eigen} we construct the simultaneous eigenvectors and in \S \ref{isom} we exhibit the isomorphism.

\subsection*{Acknowledgement}
This is a part of the master thesis written under the supervision of Professor Hiraku Nakajima. 
The author would like to thank him for his valuable comments, warm encouragement and careful proofreading.
The author also would like to thank Yoshihisa Saito for his polite answers for some questions about the q-Fock space.

\section{Preliminaries}

\subsection{Quantum toroidal algebra}\label{qta}

\subsubsection{}
In this paper we usually take $\mathcal{R}=\C(s^{1/2},t^{1/2})$ as the coefficient field.
We set
\[
p=t^l,\quad q=s^{1/2}t^{1/2},\quad r=s^{-1/2}t^{1/2}.
\]

\subsubsection{}
Let us define the \textbf{quantum toroidal algebra} $\tor$ ($l>2$).
This is an $\R$-algebra generated by $e_{i,n}$, $f_{i,n}$, $K^{\pm}_{i}$ and $h^{\pm}_{i,m}$ ($i\in I$, $n\in\Z$, $m\in\Z_{>0}$).  
The relations are expressed using the formal series
\[
e_i(z)=\sum_{n\in\Z} e_{i,n}z^{-n},\quad f_i(z)=\sum_{n\in\Z} f_{i,n}z^{-n},\\
\]
\[
K_i^{\pm}(z)=K^\pm_i\exp\left(\pm(q-q^{-1})\sum_{m>0} h^{\pm}_{i,m}z^{\mp m}\right)
\]
as follows : 
\begin{align*}
K^+_{i.0}K^-_{i.0}=K^-_{i.0}K^+_{i.0}&=1,\\
[K^\pm_i(z),K^\pm_i(w)]=[K^+_i(z),K^-_i(w)]&=0,\\
[K^\pm_i(z),e_j(w)]=[K^\pm_i(z),f_j(w)]&=0 \quad (j\neq i,i\pm 1),\\
(r^{\,\varepsilon} z-q^{-1}w)K^{\pm}_i(z)e_{i+\varepsilon}(w)&=(r^{\,\varepsilon} q^{-1}z-w)e_{i+\varepsilon}(w)K^{\pm}_i(z)\quad (\varepsilon=\pm1),\\
(z-q^{2}w)K^{\pm}_i(z)e_{i}(w)&=(q^{2}z-w)e_{i}(w)K^{\pm}_i(z),\\
(r^{\,\varepsilon} z-qw)K^{\pm}_i(z)f_{i+\varepsilon}(w)&=(r^{\,\varepsilon} qz-w)f_{i+\varepsilon}(w)K^{\pm}_i(z)\quad (\varepsilon=\pm1),\\
(q^{-2}z-w)K^{\pm}_i(z)f_{i}(w)&=(z-q^{-2}w)f_{i}(w)K^{\pm}_i(z),\\
[e_i(z),f_j(w)]&=\frac{\delta_{ij}\delta(z/w)}{q-q^{-1}}(K_i^+(w)-K_i^-(z)),\\
(r^{\,\varepsilon} z-q^{-1}w)e_i(z)e_{i+\varepsilon}(w)&=(r^{\,\varepsilon} q^{-1}z-w)e_{i+\varepsilon}(w)e_i(z)\quad (\varepsilon=\pm1),\\
(z-q^{2}w)e_i(z)e_{i}(w)&=(q^{2}z-w)e_{i}(w)e_i(z),\\
(r^{\,\varepsilon} z-qw)f_i(z)f_{i+\varepsilon}(w)&=(r^{\,\varepsilon} qz-w)f_{i+\varepsilon}(w)f_i(z)\quad (\varepsilon=\pm1),\\
(z-q^{-2}w)f_i(z)f_{i}(w)&=(q^{-2}z-w)f_{i}(w)f_i(z),\\
\{e_i(z_1)e_i(z_2)e_{i\pm1}(w)&-(q+q^{-1})e_i(z_1)e_{i\pm1}(w)e_i(z_2)\\
&+e_{i\pm1}(w)e_i(z_1)e_i(z_2)\}
+\{z_1\leftrightarrow z_2\}=0,\\
\{f_i(z_1)f_i(z_2)f_{i\pm1}(w)&-(q+q^{-1})f_i(z_1)f_{i\pm1}(w)f_i(z_2)\\
&+f_{i\pm1}(w)f_i(z_1)f_i(z_2)\}
+\{z_1\leftrightarrow z_2\}=0,
\end{align*}
where $\delta(Z)=\sum_{n\in\Z}Z^n$.

\begin{rem}
The quantum toroidal algebra in \cite{vv-k} is "twisted" in their words, which may or may not be isomorphic to ours. 
See Remark \ref{corr} for the relation between these two algebras.
\end{rem}

\subsubsection{}

The \textbf{horizontal subalgebra} $\uqc$ is the subalgebra of $\tor$ generated by $e_{i,0}$, $f_{i,0}$ and $K^{\pm}_{i}$ ($i\in I$).
This is isomorphic to $\uq\otimes\R$.

The \textbf{vertical subalgebra} $\uqd$ is the subalgebra of $\tor$ generated by $e_{i,n}$, $f_{i,n}$, $K^{\pm}_{i}$, and $h^{\pm}_{i,m}$ ($i\neq 0$, $n\in\Z$, $m\in\Z_{>0}$).
Define $\widetilde{e_{i,n}}$, $\widetilde{f_{i,n}}$, $\widetilde{K^{\pm}_{i}}$ and $\widetilde{h^{\pm}_{i,m}}$ by 
\begin{align*}
\widetilde{e_i(z)}&=\sum_{n\in\Z}\widetilde{e_{i,n}}z^{-n}=e_i(r^{-l+i}z),\\
\widetilde{f_i(z)}&=\sum_{n\in\Z}\widetilde{f_{i,n}}z^{-n}=f_i(r^{-l+i}z),\\
\widetilde{K^\pm_i(z)}&=\widetilde{K^\pm_i}\exp\left(\pm(q-q^{-1})\sum_{m>0} \widetilde{h^{\pm}_{i,m}}z^{\mp m}\right)
=K^\pm_i(r^{-l+i}z).
\end{align*}
They satisfy the relations in Drinfeld new realization of $\uq$, and so $\uqd$ is also isomorphic to $\uq\otimes\R$.

\subsection{Notations for Young diagrams}

\subsubsection{}\label{young}

Let $\Pi$ denote the set of all Young diagrams. 
We identify a Young diagram with a subset of $(\Z_{\geq0})^2$.
A {\bf node} is an element of $(\Z_{\geq 0})^2$.

The {\bf content} of a node $(x,y)$ is the number $x-y$.
A node is called an $i$-node if its content equals to $i$ modulo $l$.
For $\lambda\in\Pi$ let $d_i(\lambda)$ denote the number of $i$-nodes in $\lambda$
and set $\bal(\lambda)=(d_i(\lambda))_{i=0,\ldots,l-1}\in\Z^l$.
We define the order $>$ on the set of nodes according to their contents.

For $\lambda\in\Pi$ a node $(x,y)$ is called {\bf addable} 
if $(x,y)\notin\lambda$ and $(x-1,y),(x,y-1)\in\lambda$.
A node $(x,y)$ is called {\bf removable} 
if $(x,y)\in\lambda$ and $(x+1,y),(x,y+1)\notin\lambda$.
Let $A_{\lambda,i}$ (resp. $R_{\lambda,i}$) denote the set of all addable (removable) $i$-nodes of $\lambda$.

A \textbf{hook} is a pair $\left((x_h,y_h),(x_t,y_t)\right)$ such that 
$(x_h,y_h-1),(x_t,y_t)\in\lambda$ and
$(x_h,y_h),(x_t+1,y_t)\notin\lambda$.
The \textbf{hook length} of a hook $\left((x_h,y_h),(x_t,y_t)\right)$ is the number $-x_h+y_h+x_t-y_t$.
A hook is called an {\bf $l$-hook} if its length is a multiple of $l$.

\subsubsection{}

A \textbf{Maya diagram} with charge $c$ is 
an infinite decreasing sequence of integers $\bk=(k_1,k_2,\ldots)$ such that $k_a=-a+c$ for sufficiently large
$a$.
A Maya diagram with charge $c$ can be identified with a Young diagram
\[
\lambda=\coprod_{\begin{subarray}{c}
a\in\Z_{>0}\\
1\leq b\leq k_a+a-c
\end{subarray}}
(a-1,b-1).
\]
Let $\Pi_c$ denote the set of all Maya diagrams with charge $c$. Then $\Pi$ and $\Pi_c$ are bijective. 


\subsubsection{}

We sometimes identify a Maya diagram $\bk=(k_1,k_2,\ldots)$ with the subset $\{k_1,k_2,\ldots\}$ of $\Z$.

If $k_a-1\notin \bk$ ($a\in\Z_{>0}$), then a node $(a-1,k_a+a-c-1)$ is a removable node. 
Its content equals to $c-k_a$.
If $k_a+1\notin \bk$ ($a\in\Z_{>0}$), then a node $(a-1,k_a+a-c)$ is an addable node.
Its content equals to $c-k_a-1$.

Note that $\{(a,b)\mid a\in\bk, b\notin\bk, a>b\}$ is a finite set.
Such a pair $(a,b)$ corresponds to a hook in term of Young diagram.
Its hook length is $a-b$.

\section{K-theory of quiver varieties}\label{k-theory}
{\ }

\vspace{-14pt}

In this section we review the representation of $\tor$ on the equivariant K-groups of the quiver varieties of type $\hat{A}$.

\smallskip

A quiver variety, introduced by Nakajima, is a certain moduli space of representations of a quiver.
He also introduced a certain subvariety of the product of two quiver varieties called the Hecke correspondence.
Using the Hecke correspondence we can construct an action of the quantum affinization of the Kac-Moody algebras 
on the torus equivariant K-groups of the quiver varieties (\cite{quiver3}). 

By the localization theorem, localized equivariant K-groups have bases indexed by fixed points.
The fixed points of the quiver varieties of type $\hat{A}$ are indexed by Young diagrams.
The action of the quantum toroidal algebra can be written in terms of Young diagrams (\cite{vv-k}, see Theorem \ref{formula of toroidal action}).

\smallskip

Nakajima's definition of quiver varieties involves parameters $\dimv$ and $\dimw$, 
where $\dimw$ corresponds to the $l$-highest weight of the representation.
In this paper we work on the case $\dimw=(1,0,\ldots,0)$ only, 
in other words, we deal with the level $1$ representation only.

We do not take the original definition of quiver varieties but another equivalent one, 
which works only for the case $\dimw=(1,0,\ldots,0)$.  

We use $\bal$ instead of $\dimv$.

\subsection{Quiver varieties}\label{quivervar}

\subsubsection{}\label{qvdef}

Let $\hilb$ denote the Hilbert scheme of $n$ points on $\C^2$ :
\[
\hilb=\{J\underset{\text{ideal}}{\subset}\C[x,y]\mid \mathrm{dim}\,\C[x,y]/J=n\},
\] 
and $\mathrm{Sym}^n\C^2$ denote the $n$-th symmetric product of $\C^2$ :
\[
\mathrm{Sym}^n\C^2=\left\{\sum a_ip_i\ \big|\ a_i\in\Z_{>0},\ \sum a_i=n,\ p_i\in\C^2\right\}
\]
Let $\pi$ denote the Hilbert-Chow morphism :
\[
\begin{array}{rccc}
\pi\colon&\hilb&\longrightarrow&\mathrm{Sym}^n\C^2\\
&J&\longmapsto&\mathrm{supp}\, \C[x,y]/J.
\end{array}
\]

We regard $\Z/l\Z$ as the subgroup of $\mathrm{SU(2)}$. 
The action of $\Z/l\Z$ on $\C^2$ induces the action of $\Z/l\Z$ on $\mathrm{Sym}^n\C^2$ and $\hilb$ so that $\pi$ is $\Z/l\Z$-equivariant.
Let $\left(\mathrm{Sym}^n\C^2\right)^{\Z/l\Z}$ and $\left((\C^2)^{[n]}\right)^{\Z/l\Z}$ denote the sets of the fixed points.

Note that for $J\in\left((\C^2)^{[n]}\right)^{\Z/l\Z}$, $\C[x,y]/J$ has a canonical $\Z/l\Z$-module structure.
For $\bal=(d_0,\ldots,d_{l-1})\in\Z^l$ such that $\sum d_i=n$ we define the quiver variety $\M(\bal)$ by 
\[
\M(\bal)=\left\{J\in\left((\C^2)^{[n]}\right)^{\Z/l\Z}\ \Big |\ \mathrm{dim}\,\C[x,y]/J\simeq\bigoplus_i\C_{(i)}^{\ \oplus \alp_i}\right\},
\]
where $\C_{(i)}$ is the $1$-dimensional $\Z/l\Z$-module with weight $i$.

We set
\[
\kappa_1=\frac{1}{2}\bal C\hspace{1pt}{{}^t\hspace{-2pt}\bal}+\alp_0,\quad \kappa_2=n-\kappa_1l,
\]
where $C$ is the Cartan matrix of type $\hat{A}_{l-1}$.
Then we have $\mathrm{dim}\,\M(\bal)=2\kappa_1$ (\cite{quiver1}). 
Let $\zeta\in\Z/l\Z$ be a generator of $\Z/l\Z$.
We define the closed subvariety
\[
\M_0(\bal)=\left\{\kappa_2\,[0]+\sum_{j=1}^{\kappa_1}\left([p_j]+\cdots+[\zeta^{l-1}p_j]\right)\in\left(\mathrm{Sym}^n\C^2\right)^{\Z/l\Z}\ \Big|\ p_j\in \C^2\right\}.
\]
of $\left(\mathrm{Sym}^n\C^2\right)^{\Z/l\Z}$. Then we have $\pi(\M(\bal))\subset\M_0(\bal)$ (in fact we can check $\pi(\M(\bal))=\M_0(\bal)$).

For $\bal,\,\bal'\in\Z^l$ such that $\alp_i\leq \alp'_i$ for all $i$, we have the inclusion given by    
\[
\begin{array}{ccc}
\M_0(\bal)&\hookrightarrow&\M_0(\bal')\\
X&\longmapsto&X+\kappa_3\,[0],
\end{array}
\]
where $\kappa_3=\sum(\alp'_i-\alp_i)$.
We set 
\[
\M=\coprod_{\bal\in\Z^l}\M(\bal),\quad \M_0=\bigcup_{\bal\in\Z^l}\M_0(\bal),
\]
and
\[
\mathfrak{Z}=\M\underset{\M_0}{\times}\M.
\]
Note that we introduce $\M_0$ just only for terminological reason.
We work on $\M$ and $\mathfrak{Z}$, of which connected components are finite dimensional.

\subsubsection{}\label{torus action on qv}

The natural $T=(\C^*)^2$-action on $\C^2$ induces a $T$-action on $\M$.
The $T$-fixed points of $\M$ are indexed by $\Pi$.
For $\lambda\in\Pi$ the corresponding ideal $J_\lambda\in \left(\C^2\right)^{[\mathrm{deg}\lambda]}$ is the ideal generated by $\{x^ay^b\mid (a,b)\notin\lambda\}$.
Then $\{[x^ay^b]\in\C[x,y]/J\mid (a,b)\in\lambda\}$ forms a basis of $\C[x,y]/J$.

For $\zeta\in\Z/l\Z$ we have $\zeta\cdot[x^ay^b]=\zeta^{a-b}[x^ay^b]$.
So $J_\lambda\in\M(\bal(\lambda))$.

For $(s,t)\in T$ we have $(s,t)\cdot[x^ay^b]=s^at^b[x^ay^b]$. 
So $\C\cdot[x^ay^b]=s^at^b\in R(T)=\Z[s^\pm,t^\pm]$, 
where 
$R(T)$ is the representation ring of $T$ and 
we identify the coordinate functions of $T$ with the generators of $R(T)$. 
Thus for a node $X=(a,b)$, we set $\node{X}=s^at^b\in R(T)$.

\subsection{Representation on K-theory of quiver varieties}\label{toroidal on K}

\subsubsection{}\label{Hecke}

Let ${\bf e}_i$ denote the $i$-th coordinate vector in $\Z^l$.
For $\bal\in \Z^l$ we define the  subvariety of $\mathfrak{Z}$ by 
\[
\B_i(\bal)=\left\{(J_1,J_2)\in\qZ\mid J_1\in\M(\bal),\ J_2\in\M(\bal+{\bf e}_i),\ J_1\supset J_2\right\}
\]
This is called the {\bf Hecke correspondence}.

Let $p_\varepsilon$ denote the projection from $\qZ$ to the $\varepsilon$-th factor ($\varepsilon=1,2$) 
and $q_\varepsilon$ denote its restriction to
$\B_i(\bal)\subset\qZ$. 
We define the tautological bundle $\tL$ on $\B_i(\bal)$ by $q_2^*\mathfrak{V}/q_1^*\mathfrak{V}$.

\subsubsection{}\label{corr}

For a $T$-equivariant vector bundle $\mathfrak{B}$ on $X$,
let $\mathrm{det}\mathfrak{B}$ denote its determinant, 
$\wedge^i\mathfrak{B}$ denote its $i$-th wedge product, 
and set $\wedge_z\mathfrak{B}=\sum_{i\geq 0}(-z)^i\wedge^i\mathfrak{B}$.
These operators can be extended to operators on $\mathrm{K}^T(X)$.
For a $\Z/l\Z$-module $M$ we set $M_i=\mathrm{Hom}_{\Z/l\Z}(\C_{(i)},M)$.

We set 
\[
\mathfrak{H}=(-1+s+t-st)\mathfrak{V}+\mathfrak{W}\in \mathrm{K}^T(\M).
\] 
We define an action of $\tor$ on $\mathrm{K}^T_\R(\M)=\mathrm{K}^T(\M)\otimes \R$ by
\begin{align*}
e_{i,n}(x)&=c_i^-(\bal)\,{p_1}_{*}\left(p_2^{\ *}x\otimes(\tL)^{n+h_i(\bal)}\right)\quad&x\in \mathrm{K}^T_\R(\M(\bal+{\bf e}_i)),\\
f_{i,n}(x)&=c_i^+(\bal)\,{p_2}_{*}\left(p_1^{\ *}(x\otimes \mathrm{det}(s^{-1}t^{-1}\mathfrak{H}_i)\otimes\tL^{n}\right)\quad&x\in \mathrm{K}^T_\R(\M(\bal)),\\
K^\pm_i(z)(x)&=c_i^-(\bal)c_i^+(\bal)\left(\wedge_z\hspace{-2pt}\left((s^{-1}t^{-1}-1)\mathfrak{H}_i^*\right)\right)^\pm x\quad&x\in \mathrm{K}^T_\R(\M(\bal)),
\end{align*}
where
the index ${}^+$ (resp. ${}^-$) means the expansion as a formal power series in $z^{-1}$ (resp. $z$) and
\begin{align*}
c_i^-(\bal)&=(-1)^{\alp_i}\,s^{\,(2\alp_i-\alp_{i+1}+1)/2}\,t^{\,(-2\alp_{i-1}+2\alp_i-\alp_{i+1}+1)/2},\\
c_i^+(\bal)&=(-1)^{-\alp_{i-1}+\alp_{i}-\alp_{i+1}}\,s^{-\alp_{i-1}/2}\,t^{\,\alp_{i-1}/2},\\
h_i(\bal)&=\alp_{i-1}-2\alp_i+\alp_{i+1}.
\end{align*}

\begin{rem}
We slightly modify the actions in \cite{vv-k}. In fact we have
\begin{align*}
e_{i,n}&=(-1)^{\alp_{i+1}}s^{(\alp_{i+1}+1)/2}t^{(-\alp_{i+1}+1)/2}\Omega^-_{i,n},\\
f_{i,n}&=(-1)^{\alp_{i+1}}s^{-\alp_{i+1}/2}t^{\alp_{i-1}/2}\Omega^+_{i,n},\\
K^\pm_i(z)&=s^{(-\alp_{i-1}+\alp_{i+1}+1)/2}t^{(\alp_{i-1}-\alp_{i+1}+1)/2}\Theta^\pm_{i}(z).
\end{align*}
Here the operators on right hand side are defined in 3.3 of \cite{vv-k}, 
where we should replace their symbols $q$, $t$, $k$, $s$ with our symbols $t$, $s$, $i$, $n$.

Substitute this to theorem 2 in \cite{vv-k} and the definition of quantum toroidal algebra in \cite{vv-k},
we can verify $e_{i,n}$, $f_{i,n}$ and $K^\pm_i(z)$ satisfy the relation in \ref{qta}.
\end{rem}

\subsubsection{}\label{formula of toroidal action}

Let $i_{\lambda}$ denote the inclusion $\{J_\lambda\}\hookrightarrow\M$
and $1_{\lambda}$ denote the generator of $\mathrm{K}^{T}(\{J_\lambda\})$.
We set 
$b_\lambda={i_{\lambda}}_*(1_{\lambda})\in \mathrm{K}^{T}(\M)$.

By the localization theorem 
\[
\mathrm{K}^T_\mathcal{R}(\M)\simeq \bigoplus_{\lambda\in\Pi}\mathcal{R}b_\lambda.
\]

\begin{thm}(\cite{vv-k} lemma 8)
For $\lambda\in\Pi$ such that $\mathbf{d}(\lambda)=\mathbf{d}$ we have
\begin{align*}
e_{i,n}(b_{\lambda})=&\,
(-s^{1/2}t^{-1/2})^{\alp_{i-1}}\sum_{X\in R_{\lambda,i}}\hspace{2pt}
\Biggl[\node{X}^{n}\hspace{-2pt}
\prod_{A\in A_{\lambda,i}}\left((st)^{1/2}\node{A}^{*}-(st)^{1/2}\node{X}^{*}\right)^{-1}
\\&
\quad\quad\quad\quad\quad\quad\quad\quad\quad
\times\prod_{R\in R_{\lambda\backslash X,i}}\left((st)^{-1/2}\node{R}^*-(st)^{1/2}\node{X}^{*}\right)\,b_{\lambda\backslash X}\Biggr],\\
f_{i,n}(b_{\lambda})=&\,
(-s^{1/2}t^{-1/2})^{-\alp_{i-1}}\sum_{X\in A_{\lambda,i}}\hspace{2pt}\Biggl[\node{X}^{n}\hspace{-2pt}
\prod_{A\in A_{\lambda\cup X,i}}\left(st\node{A}^{*}-\node{X}^{*}\right)
\\&
\quad\quad\quad\quad\quad\quad\quad\quad\quad\quad\quad\quad\quad
\times\prod_{R\in R_{\lambda,i}}\left(\node{R}^*-\node{X}^{*}\right)^{-1}b_{\lambda\cup X}\Biggr],\\
K_i^{\pm}(z)(b_{\lambda})=&
\left(\prod_{A\in A_{\lambda,i}}\frac{(st)^{1/2}\node{A}^{*}z-(st)^{-1/2}}{\node{A}^{*}z-1}
\prod_{R\in R_{\lambda,i}}\frac{(st)^{-1/2}\node{R}^{*}z-(st)^{1/2}}{\node{R}^{*}z-1}\right)^{\pm}b_{\lambda}.
\end{align*}
where $\left(s^at^b\right)^*=s^{-a}t^{-b}$ for $s^at^b\in R(T)$.
\end{thm}

\section{Schur-Weyl duality}\label{schur}

{\ }

\vspace{-14pt}

In this section we review Schur-Weyl duality.

One can construct representations of the quantum affine algebra $\uq$ from representations of the affine Hecke algebra $\aha$ (\cite{chari-pressley-schur-weyl}, see \ref{CP}).
In this construction the action is given originally in terms of Chevalley generators.
One can rewrite the action in terms of Drinfeld generators (\cite{vv-schur}, see Theorem \ref{DCP}).

Further, Schur-Weyl duality in \cite{chari-pressley-schur-weyl} can be extended to get representations of the quantum toroidal algebra $\tor$ from representations of the toroidal Hecke algebra $\daha$.
This is done by extending the action of $\uq$ to $\tor$ using the rotation automorphism of the Dynkin diagram of type $\hat{A}$ ([VV1], see Theorem \ref{VVschur}).

\subsection{Schur-Weyl duality for affine algebras}\label{schur for affine}

\subsubsection{}\label{affine hecke}

The \textbf{finite Hecke algebra} $\ha$ is the $\mathcal{R}$-algebra generated by 
${T_a}^{\pm1}\ (a=1,\ldots,N-1)$ with relations : 
\begin{align*}
T_a{T_a}^{-1}={T_a}^{-1}T_a&=1,\\
(T_a+1)(T_a-q^2)&=0,\\
T_aT_{a+1}T_a&=T_{a+1}T_aT_{a+1},\\
T_aT_b&=T_bT_a\quad (|a-b|>1).
\end{align*}
The \textbf{affine Hecke algebra} $\aha$ is the $\mathcal{R}$-algebra generated by 
${T_a}^{\pm1}\ (a=1,\ldots,N-1)$, ${X_a}^{\pm1}\ (a=1,\ldots,N)$ with relations : 
\begin{align*}
T_a{T_a}^{-1}={T_a}^{-1}T_a&=1,\\
(T_a+1)(T_a-q^2)&=0,\\
T_aT_{a+1}T_a&=T_{a+1}T_aT_{a+1},\\
T_aT_b&=T_bT_a\quad (|a-b|>1),\\
X_aX_b&=X_bX_a,\\
T_aX_aT_a&=q^2X_{a+1},\\
X_bT_a&=T_aX_b\quad (b\neq a,a+1).
\end{align*}

\subsubsection{}\label{left action}

Let $V=\mathcal{R}^l$ with a basis $\{v_0,\ldots,v_{l-1}\}$.
We define $\vt\in\mathrm{End}(V^{\otimes  2})$ by
\[
\vt(v_{i_1}\otimes v_{i_2})=
\begin{cases}
q^2v_{i_1}\otimes v_{i_2} & \text{if $i_1=i_2$},\\
qv_{i_2}\otimes v_{i_1} & \text{if $i_1<i_2$},\\
qv_{i_2}\otimes v_{i_1}+(q^2-1)v_{i_1}\otimes v_{i_2} & \text{if $i_1>i_2$}.
\end{cases}
\]
Then we have a left action of $\ha$ on $V^{\otimes N}$ defined by 
\[
T_a\longmapsto \vt_a=1^{\otimes a-1}\otimes \vt\otimes 1^{\otimes N-a-1}.
\]

\subsubsection{}\label{CP}

Let $M$ be a right $\aha$-module. 
We define the following operators on $M\otimes_{\ha}V^{\otimes N}$ : 
\begin{align*} 
e_i(m\otimes v)&=\sum_{a=1}^N mX_a^{\delta_{i,0}}\otimes (K^i_1)^{-1}\cdots (K^i_{a-1})^{-1}E_a^{i,i-1}\,v,\\ 
f_i(m\otimes v)&=\sum_{a=1}^N mX_a^{-\delta_{i,0}}\otimes E_a^{i-1,i}K^i_{a+1}\cdots K^i_N\,v,\\
h_i(m\otimes v)&=m\otimes K_1^i\cdots K_N^i\,v.
\end{align*} 
Here $ E_a^{i,j} = 1^{\otimes^{a-1}}\otimes E^{i,j} \otimes 1^{\otimes^{N-a}},$ 
where $E^{i,j}\in\mathrm{End}(V)$ is the matrix unit with respect to the basis $v_0,\ldots,v_{l-1}$ and 
$K_a^i = q^{E_a^{i-1,i-1} - E_a^{i,i}}$. 
These operators give a left $\uq$-action on $M\otimes_{\ha}V^{\otimes N}$ (\cite{chari-pressley-schur-weyl}).

\subsubsection{}\label{DCP}

An isomorphism between the algebras defined by Chevalley generators
and by Drinfeld new realization is given in \cite{beck-newrealization}.

For $\bj=(j_1,\ldots,j_N)\in\{0,\ldots,l-1\}^N$ let $\bv_\bj$ denote $v_{j_1}\otimes\cdots\otimes v_{j_N}\in V^{\otimes N}$.

For $1\leq a,b \leq N$ we define
\[
T_{a,b}=
\begin{cases}
T_{a}T_{a+1}\cdots T_{b-1} & a<b,\\
1 & a=b,\\
T_{a-1}T_{a-2}\cdots T_{b} & a>b.
\end{cases}
\]

\begin{thm}\label{lemma of VV}
(\cite{vv-schur} Theorem 3.3)
Assume $\bj$ is an non-decreasing sequence.
We put $n_i=\sharp\{a\mid j_a=i\}$ and $\bar{n}_i=\sum_{i'=1}^in_{i'}$.  Let us write $\bj=[n_0,n_1,\ldots]$.

For $m\otimes\bv_\bj\in M\otimes_{\ha}V^{\otimes N}$ the actions of Drinfeld generators of $\uq$ are described as follows :
\begin{align*}
\widetilde{e_i(z)}(m\otimes\bv_\bj)&=q^{1-n_{i}}m\left(\sum_{a=\bar{n}_{i-1}+1}^{\bar{n}_{i}}T_{a,\bar{n}_{i-1}+1}\right)\delta\left(q^{l-i}Y_{\bar{n}_{i-1}+1}z\right)\otimes\bv_{\bj^-},\\
\widetilde{f_i(z)}(m\otimes\bv_\bj)&=q^{1-n_{i-1}}m\left(\sum_{a=\bar{n}_{i-2}+1}^{\bar{n}_{i-1}}T_{a,\bar{n}_{i-1}}\right)\delta\left(q^{l-i}Y_{\bar{n}_{i-1}}z\right)\otimes\bv_{\bj^+},\\
\widetilde{K^{\pm}_i(z)}(m\otimes\bv_\bj)&=m\prod_{j_a=i-1}\theta_1^{\pm}\left(q^{l-i+1}Y_{a}z\right)\prod_{j_b=i}\theta_{-1}^{\pm}\left(q^{l-i-1}Y_{b}z\right)\otimes\bv_\bj,
\end{align*}
where $\bj_-=[\ldots,n_{i-1}+1,n_{i}-1,\ldots]$, $\bj_+=[\ldots,n_{i-1}-1,n_{i}+1,\ldots]$ and $\theta_m(z)=\frac{q^mz-1}{z-q^m}$.
\end{thm}

\subsection{Schur-Weyl duality for toroidal algebras}\label{schur for toroidal}

\subsubsection{}\label{toroidal hecke}

The \textbf{toroidal Hecke algebra} $\daha$ is the $\mathcal{R}$-algebra generated by
${T_a}^{\pm1}\ (a=1,\ldots,N-1)$, ${X_a}^{\pm1}\ (a=1,\ldots,N)$, ${Y_a}^{\pm1}\ (a=1,\ldots,N)$ 
with relations : 
\begin{align*}
T_a{T_a}^{-1}={T_a}^{-1}T_a&=1,\\
(T_a+1)(T_a-q^2)&=0,\\
T_aT_{a+1}T_a&=T_{a+1}T_aT_{a+1},\\
T_aT_b&=T_bT_a\quad (|a-b|>1),\\
X_aX_b&=X_bX_a,\\
T_aX_aT_a&=q^2X_{a+1},\\
X_bT_a&=T_aX_b\quad (b\neq a,a+1),\\
Y_aY_b&=Y_bY_a,\\
T_a^{-1}Y_aT_a^{-1}&=q^{-2}Y_{a+1},\\
Y_bT_a&=T_aY_b\quad (b\neq a,a+1),\\
X_0Y_1&=pY_1X_0,\\
X_2Y_1^{-1}X_2^{-1}Y_1&=q^{-2}T_1^2,
\end{align*}
where $X_0=X_1\cdots X_N$.

Let $\ahay$ (resp. $\ahax$) denote the subalgebra generated by $\{{T_a}^{\pm 1}\}$ and $\{Y_a\}$
(resp. $\{{T_a}^{\pm1}\}$ and $\{X_a\}$). They are isomorphic to $\aha$.

\subsubsection{}\label{VVschur}

Let $M$ be a right $\daha$-module.
Regarding $M$ as a right $\ahay$-module
we have the action of $\uq$ on $M\otimes_{\ha}V^{\otimes N}$ by \ref{CP}.

We define an operator $\rho$ on $M\otimes_{\ha}V^{\otimes N}$ by
\[
\rho(m\otimes v_{i_1}\otimes\cdots\otimes v_{i_N})
=mX_1^{\delta_{0,l_1}}\cdots X_N^{\delta_{0,i_N}}\otimes v_{i_1-1}\otimes\cdots\otimes v_{i_N-1}.
\]

\begin{lem}(\cite{vv-schur} Proposition 3.4)
We set
$\mathcal{X}_i(z)=\widetilde{\mathcal{X}_i(r^{\,l-i}z)}$ ($\mathcal{X}=e,f,K^{\pm}$). Then we have
\[
\mathcal{X}_{i-1}(z)={\rho}^{-1}\circ\mathcal{X}_{i}(q^{-1}r^{-1}z)\circ{\rho}
\]
\end{lem}

\begin{thm}(\cite{vv-schur} Theorem 3.5)
The action of $\uq\otimes\R\simeq \uqd$ on $M\otimes_{\ha}V^{\otimes N}$ can be extended to an action of 
$\tor$ so that the actions of
$\mathcal{X}_0$ ($\mathcal{X}=e,f,K^{\pm}$) are   
are given by
\[
\mathcal{X}_0 (z)={\rho}^{-1}\circ\mathcal{X}_{1}(q^{-1}r^{-1}z)\circ{\rho}
\]
\end{thm}




\section{Representation on the q-Fock space}\label{toroidal on q}

{\ }

\vspace{-14pt}

In this section we review the action of $\tor$ on the q-Fock space following \cite{saito-takemura-uglov} and \cite{vv-q-fock}.

As a q-analogue of the permutation representation ,
$\mathcal{R}[z_1^\pm,\ldots,z_N^\pm]$ has a right $\ha$-module structure.
We define the q-wedge space by $\mathcal{R}[z_1^\pm,\ldots,z_N^\pm]\otimes_{\ha} V^{\otimes N}$.
This is the q-analogue of the classical wedge space $\otimes^N V(z)/\oplus\mathrm{Ker}(\mathrm{id}+\sigma_i)$, 
where $\sigma_i$ is the generator of $\mathfrak{S}_N$.

We define the q-Fock space taking "limit" of the q-wedge space. 
In other words the q-Fock space is the q-analogue of the classical semi-infinite wedge space.  

It is known
the right $\ha$-module structure on $\mathcal{R}[z_1^\pm,\ldots,z_N^\pm]$ can be extended to a right $\daha$-module structure called Dunkl-Cherednik representation.
By Schur-Weyl duality described in \ref{VVschur}, we have an action of $\tor$ on the q-wedge space.
This can be naturally lifted to an action on the q-Fock space.

\subsection{The q-Fock space}\label{q-wedge and q-fock}

Here we review the definition of the q-Fock space. The reader can refer to \cite{kashiwara-miwa-stern} for detail.

\subsubsection{}\label{right action}

For $1\leq a<b\leq N$
let us define an operator $g_{ab}$ on 
$\mathcal{R}[z_1^{\pm 1},\ldots ,z_N^{\pm 1}]$ by  
\[
g_{ab}=\frac{q^{-1}z_a-qz_b}{z_a-z_b}(\sigma_{ab}-1)+q,
\]
where $\sigma_{ab}$ is the operator defined by the permutation of variables $z_a$ and $z_b$.

Then we have a right action of $\ha$ on $\poly$ defined by
\[
T_a\longmapsto \pt_a=(q^2-1)-qg_{a,a+1}.
\]

\subsubsection{}\label{q-wedge}

Let $V(z)=\mathcal{R}[z^{\pm 1}]\otimes V$.
We define 
\begin{align*}
\wedge ^NV(z)&=\poly \otimes_{\ha} V^{\otimes N}\\
&=\otimes^NV(z)\Big/\sum_{a=1}^{N-1}\mathrm{Im}\left(\pt_a\otimes 1_{V^{\otimes N}}-1_{\poly}\otimes\vt_a\right).
\end{align*}
This is called the \textbf{q-wedge space}.

\subsubsection{}\label{q-wedge relation}

We write $u_k=z^m\otimes v_j$ for $k=j-l(m+1)$.
Let $u_{k_1}\wedge\cdots\wedge u_{k_N}$ denote the image of $u_{k_1}\otimes\cdots\otimes u_{k_N}$ for the quotient map.
We say $u_{k_1}\wedge\cdots\wedge u_{k_N}$ is \textbf{normally ordered} if $k_a>k_b$ for $a<b$.

For $N=2$ we can verify that if $k\equiv k'$ then 
\[
u_k\wedge u_{k'}=-u_{k'}\wedge u_k,
\]
and if $k<k'$ and $k-k'\equiv i\ (1\leq i\leq l-1)$ then 
\begin{align*}
u_k\wedge u_{k'}=-qu_{k'}\wedge u_k+&(q^2-1)\bigl(u_{k'-i}\wedge u_{k+i}\\
&-qu_{k'-l}\wedge u_{k+l}+q^2u_{k'-l-i}\wedge u_{k+l+i}-\cdots\bigr)
\end{align*}
where the summation continues as long as the wedge is normally ordered. 

The set of all normally ordered wedges forms a basis of $\wedge ^NV(z)$.





\subsubsection{}\label{iota}

For $c\in\Z$ and $0<N<N'$ we define
\[
\begin{array}{rccc}
\iota^c_{N,N'}\colon&\Z^N&\longrightarrow&\Z^{N'}\\
&(k_1,\ldots,k_N)&\longmapsto&(k_1,\ldots,k_N,-N+c-1,\ldots,-N'+c).
\end{array}
\]
For $\bk=(k_1,\ldots,k_N)\in\Z^N$ let us write $u_\bk=u_{k_1}\wedge\cdots\wedge u_{k_N}$.
We can check the well-definedness of the map
\[
\begin{array}{ccl}
\wedge^NV(z) & \longrightarrow & \wedge^{N'}V(z)\\
u_\bk & \longmapsto & u_{\iota^c_{N,N'}(\bk)}.
\end{array}
\]
We write $\iota^c_{N,N'}$ for this map as well.

\subsubsection{}\label{q-fock}

We define 
\[
F(c)=\varinjlim_{\iota^c_{N,N'}}\wedge^NV(z),\quad F=\bigoplus_{c\in\Z}F(c),
\]
and $\iota^c_{N,\infty}$ by the canonical map from $\wedge^NV(z)$ to $F(c)$.
$F$ (resp. $F(c)$) is called the \textbf{q-Fock space}  (with charge $c$).
An element of $F$ (resp. $F(c)$) is called a \textbf{semi-infinite wedge} (with charge $c$).

Let $\bk=(k_1,k_2,\ldots)$ be a Maya diagram with charge $c$ 
(we use $\bk$ both for an element of $\Z^N$ and for an infinite sequence of integers by abuse of notations), 
then $u_\bk=u_{k_1}\wedge u_{k_2}\wedge\cdots$ is a semi-infinite wedge with charge $c$. 
Note that $\{u_\bk\mid \bk\in\Pi_c\}$ forms a basis of $F(c)$.

\subsection{Representation on the q-Fock space}\label{toroidal on q-fock}

\subsubsection{}\label{dunkel rep}

Let us consider the following operators on $\mathcal{R}[z_1^\pm,\ldots,z_N^\pm]$ :  
\[
Y_a^{(N)} = g_{a,a+1}^{-1}\sigma_{a,a+1} \cdots g_{a,N}^{-1} \sigma_{a,N} p^{D_a} \sigma_{1,a} g_{1,a} \cdots \sigma_{a-1,a}g_{a-1,a}\quad (a\in\{1,\ldots,N\})
\]
where $p^{D_a}$ is the difference operator given by 
\[
p^{D_a}f(z_1,\ldots,z_a,\ldots,z_N)  =  f(z_1,\ldots,pz_a,\ldots,z_N) , \quad f \in \mathcal{R}[z_1^{\pm 1},\ldots, z_N^{\pm 1}].  
\]
The operator $Y_a^{(N)}$ is called \textbf{Dunkl-Cherednik operator}.
Then the action of $\aha$ defined in \ref{left action} can be extended to the action of $\daha$ by
\[
T_a\longmapsto \pt_a,\quad X_i\longmapsto z_a, \quad Y_a\longmapsto q^{1-N}Y^{(N)}_a.
\]    
This is called \textbf{Dunkl-Cherednik representation} (\cite{cherednik-dunkl}, \cite{cherednik-macdonald-conj}, \cite{cherednik-nsmp}).

By the Schur-Weyl duality explained in \ref{VVschur},
we have an action of $\tor$ on $\wedge ^NV(z)=\poly \otimes_{\ha} V^{\otimes N}$.

\subsubsection{}\label{TU}

For $\bk\in\Z^N $ we define $\bm\in\Z^N$ and $\bj\in\{0,\ldots,l-1\}$  by $k_a=j_a-l(m_a+1)$. 
Note that $z^{\bm}\otimes \bv_{\bj}=u_\bk$.
We identify $\bk\in\Z^N$ with the pair $(\bm,\bj)$.
Let $\bm^c=(m_1,\ldots,m_N)\in\Z^N$ denote the sequence obtained from $\bk^c=(c-a)_{1\leq a\leq N}$.

Let $\mathcal{M}^{\,c,r}_{N,l}$ denote the set of all $\bm$ such that 
\begin{itemize}
\item $\bm$ is non-decreasing with no more than $l$ elements of any given value, and  
\item $m_a\geq m_a^c$ for all $a$ and $\sum(m_a-m_a^c)=\gamma$.
\end{itemize}
For $\bm\in\mathcal{M}^{\,c,\gamma}_{N,l}$ we define
\[
\mathcal{J}(\bm)=\{\bj\in\{0,\ldots,l-1\}^N \mid \text{$j_a<j_b$ for $a<b$ such that $m_a=m_b$}\}.
\]

We define
\[
V_N^{c,\gamma}=\bigoplus_{\bm\in\mathcal{M}^{\,c,\gamma}_{N,l}}\ \bigoplus_{\bj\in\mathcal{J}(\bm)}\mathcal{R} u_{\bk}\subset\wedge^NV(z).
\]
We can check this is invariant under the $\uqd$-action.

We can see that for $\alpha,\,\beta\in\Z$ such that $\alpha l+c>\gamma l$ and $\beta>\alpha$ the restriction
\[
\iota^{c,\gamma}_{\alpha l+c,\beta l+c}=\iota^{c}_{\alpha l+c,\beta l+c}|_{V_{\alpha l+c}^{c,\gamma}}\colon V_{\alpha l+c}^{c,\gamma}\longrightarrow V_{\beta l+c}^{c,\gamma}
\]
is an isomorphism as vector spaces.

\begin{thm}(\cite{takemura-uglov-level0} Proposition 6)
$\iota^{c,\gamma}_{\alpha l+c,\beta l+c}$ is an isomorphism as $\uqd$-modules.
\end{thm}

\subsubsection{}\label{corTU}

For $\maya\in\Pi^c$ we set $\mathrm{deg}\,\maya=\sum(m_a-m_a^c)$. Note that this is well-defined. 
We set
\[
F(c)_\gamma =\bigoplus_{
\begin{subarray}{c}
\gamma\in\Pi^c\\
\mathrm{deg}\,\bk=\gamma
\end{subarray}
}\R u_\bk\subset F(c).
\]
For $\alpha\in\Z$ such that $\alpha l+c>\gamma l$ the restriction
\[
\iota^{c,\gamma}_{\alpha l+c,\infty}=\iota^{c}_{\alpha l+c,\infty}|_{V_{\alpha l+c}^{c,\gamma}}\colon V_{\alpha l+c}^{c,\gamma}\longrightarrow F(c)_\gamma \\
\]
is an isomorphism as vector space. 
By Theorem \ref{TU} we can extend the $\uqd$-action to $F(c)_r$, and so to $F(c)$.

\subsubsection{}\label{toroidal on q-fock sp}

We define 
\[
\begin{array}{rccc}
\rho_N\colon & \Z^N & \longrightarrow & \Z^N\\
& (k_1,\ldots,k_N) & \longmapsto & (k_1-1,\ldots,k_N-1).
\end{array}
\]
We write $\rho_N$ as well for the map $\wedge^NV(z)\to\wedge^NV(z)$ given by $u_\bk\mapsto u_{\rho_N(\bk)}$.
We can see this is compatible with the construction of $\rho$ in \ref{VVschur}.

We also define
\[
\begin{array}{rccc}
\rho_\infty\colon & \Pi_c & \longrightarrow & \Pi_{c-1}\\
& (k_1,k_2,\ldots) & \longmapsto & (k_1-1,k_2-1,\ldots).
\end{array}
\]
and $\rho_\infty\colon F(c)\to F(c-1)$.

For $0<N<N'\leq\infty$ we have
\[
\iota^{c-1}_{N,N'}\circ\rho_{N}=\rho_{N'}\circ\iota^c_{N,N'}.
\]
Thus the action of $\tor$ on $\wedge^NV(z)$ can be extended to $F$ 
so that 
\[
X_0(z)=\rho_{\infty}^{-1}\circ X_1(q^{-1}r^{-1}z)\circ\rho_\infty\quad (X=e,f,K^{\pm}).
\]

\section{Simultaneous eigenvectors}\label{eigen}

{\ }

\vspace{-14pt}

In this section we construct simultaneous eigenvectors for the actions of $K^\pm_{i}(z)$'s on the q-Fock space, which are the main subjects of this paper.

In \S \ref{NSMP} we review nonsymmetric Macdonald polynomials (\cite{cherednik-nsmp}, \cite{macdonald-bourbaki}, \cite{opdam-nsmp}).
For $\bm\in\Z^N$, nonsymmetric Macdonald polynomial $\Phi^\bm\in\C[z_1^\pm,\ldots,z_N^\pm]$ is a simultaneous eigenvector of Dunkl-Cherednik operators.
The transition matrix between monomials and nonsymmetric Macdonald polynomials is upper triangular with respect to the Bruhat order on $\Z^N$.
The actions of the finite Hecke algebra generators $T_a$ on nonsymmetric Macdonald polynomials can be simply described (see Proposition \ref{prop of NSMP}).

For $\bk=(k_1>\cdots >k_N)$ we define a vector $\Psi^\bk=\Phi^{\ubm}\otimes \bv_{\ubj}$ in the q-wedge space (Definition \ref{eigenvec}),
where $\ubm$ and $\ubj$ are given by "renumbering" of $(k_1,\ldots,k_N)$ so that $\ubj$ is non-decreasing (\ref{renumbering}).
It follows immediately from Theorem \ref{DCP} that $\Psi^\bk$ is a simultaneous eigenvector for the actions of $K^\pm_{i}(z)$'s ($i\neq 0$). 
We can check
\begin{itemize}
\item the eigenvalues are multiplicity free (Proposition \ref{mult-free}), and
\item the transition matrix between normally ordered wedges and $\{\Psi^\bk\}$ is upper triangular (Proposition \ref{prop2}), in particular $\{\Psi^\bk\}$ forms a basis of the q-wedge space.
\end{itemize}
So the vector $\Psi^\bk$ is characterized in term of the actions of $K^\pm_{i}(z)$'s ($i\neq 0$) (Corollary \ref{mult-free}).
Further, using them we can verify 
\begin{itemize}
\item $\Psi^\bk$ is also a simultaneous eigenvector for $K^\pm_{0}(z)$ (Corollary \ref{prop3}), and
\item $\Psi^\bk$ can be lifted to the q-Fock space (Definition \ref{lift}).
\end{itemize}
We can see the eigenvalues coincide with the eigenvalues of the torus fixed points in the representation on the equivariant K-groups of the quiver varieties.

\subsection{Nonsymmetric Macdonald polynomials}\label{NSMP}

\subsubsection{}\label{order}

Let us define the \textbf{Bruhat order}. This is the partial order on $\Z^N$ given by the transitive closure of the following two relations :

For ${\bf x}=(x_1,\ldots,x_N)\in\Z^N$
\begin{itemize}
\item If $1\leq i<j\leq N$ and $x_i>x_j$ then ${\bf x}\succ\sigma_{ij}{\bf x}$, and  
\item $1\leq i<j\leq N$ and $x_i-x_j>1$ then $\sigma_{ij}{\bf x}\succ{\bf x}+{\bf e}_i-{\bf e}_j$ where ${\bf e}_i$ is the $i$-th coordinate vector. 
\end{itemize}

\subsubsection{}\label{sigma}

For ${\bf x}\in\Z^N$
let $\sigma_{\bf x}$ denote the unique element of $\mathfrak{S}_N$ satisfying the following conditions : 
\begin{itemize}
\item if $\sigma_{\bf x}(a)<\sigma_{\bf x}(a')$ then $x_{\sigma_{\bf x}(a)}\geq x_{\sigma_{\bf x}(a')}$, and
\item if $a<a'$ and $x_{\sigma_{\bf x}(a)}=x_{\sigma_{\bf x}(a')}$ then $\sigma_{\bf x}(a)<\sigma_{\bf x}(a')$.
\end{itemize}

\subsubsection{}\label{def of NSMP}

We can see that $\poly$ has the basis $\{\Phi^\bm(z)\mid\bm\in\Z^N\}$ such that
\begin{itemize}
\item $\Phi^\bm(z)=z^\bm+\sum_{\bn\prec\bm}c(\bm,\bn)z^\bn\quad (\exists c(\bm,\bn)\in\mathcal{R})$,
\item $\Phi^\bm(z)Y^{(N)}_a=\zeta_a(\bm)\Phi^\bm(z)$, where $\zeta_a(\bm)=p^{m_a}q^{2\sigma_\bm(a)-N-1}$.
\end{itemize}
$\Phi^\bm(z)$ is called \textbf{nonsymmetric Macdonald polynomial} (\cite{cherednik-nsmp}, \cite{macdonald-bourbaki}, \cite{opdam-nsmp}).



\subsubsection{}\label{prop of NSMP}

\begin{prop}(see \cite{uglov-hep-th} \S 1.5)
\[
\Phi^\bm(z)\cdot\pt_a=
\begin{cases}
\frac{(-q^2+1)}{x-1}\Phi^\bm(z)-\frac{(x-q^2)(q^2x-1)}{(x-1)^2}\Phi^{\sigma_a\bm}(z) & (m_a>m_{a+1}),\\
\frac{(-q^2+1)}{x-1}\Phi^\bm(z)                  & (m_a=m_{a+1}),\\
\frac{(-q^2+1)}{x-1}\Phi^\bm(z)-\Phi^{\sigma_a\bm}(z) & (m_a<m_{a+1}),
\end{cases}
\]
where $\sigma_a\bm=(\ldots,m_{a+1},m_a,\ldots)$ and $x=\frac{\zeta_{a+1}(\bm)}{\zeta_a(\bm)}$.
\end{prop}

\subsection{Simultaneous eigenvectors and its properties}\label{eigen and prop}

\subsubsection{}\label{renumbering}

For $\sigma\in\mathfrak{S}_N$ we define $\hat{\sigma}\in\mathfrak{S}_N$ by $\hat{\sigma}(a)=N-\sigma(a)+1$.

For $\bk\in\Z^N_+=\{\bk\in\Z^N\mid k_1>\cdots>k_N\}$
we define 
$\ubm=(\um{1},\ldots,\um{N})\in\Z^N$ and $\ubj=(\uj{1},\ldots,\uj{N})\in\{0,\ldots,l-1\}^N$ by 
\[
\um{a}=m_{\sj(a)},\quad \uj{a}=j_{\sj(a)}.
\] 
Note that
\begin{itemize}
\item $\bm$ is non-decreasing, and if $a<b$, $m_a=m_b$ then $j_a>j_b$,
\item $\ubj$ is non-decreasing, and if $a<b$, $\uj{a}=\uj{b}$ then $\um{a}>\um{b}$,
\item ${\sj}^{\,-1}=\sm$.
\end{itemize}

\begin{exa}
For $\bk=(5,3,1,-6,-7,-8,-9,-10)$ we have

\begin{align*}
\left(
\begin{array}{c}
\bm\\
\bj
\end{array}
\right)
&=
\left(
\begin{array}{cccccccc}
-2&-1&-1&1&1&1&1&1\\
0&3&1&4&3&2&1&0
\end{array}
\right),\\
\left(
\begin{array}{c}
\ubm\\
\ubj
\end{array}
\right)
&=
\left(
\begin{array}{cccccccc}
1&-2&1&-1&1&1&-1&1\\
0&0&1&1&2&3&3&4
\end{array}
\right).
\end{align*}
In the following figure,
\begin{itemize}
\item enumerate the boxes from lower rows to upper rows and from right to left in a row, 
then $(j_a,m_a)$ is the coordinate of the $a$-th box, and
\item enumerate the boxes from left columns to right columns and from the top to the bottom in a column, 
then $(\uj{a},\um{a})$ is the coordinate of the $a$-th box.
\end{itemize}
\[
\begin{array}{ccccccccc}
 \vdots   &\  &\vdots&\vdots&\vdots&\vdots&\vdots&   &   
\vspace{1pt}\\
    1     & & \fbox{$-10$} &  \fbox{$-9$}  &  \fbox{$-8$}  &  \fbox{$-7$}  &  \fbox{$-6$}  &   &   
\vspace{2pt}\\
    0     & &  -5  &  -4  &  -3  &  -2  &  -1  &   &   
\vspace{2pt}\\    
   -1     & &   0  &   \fbox{$1$}  &   2  &   \fbox{$3$}  &   4  &   &   
\vspace{2pt}\\
   -2     & &   \fbox{$5$}  &   6  &   7  &   8  &   9  &   &   
\vspace{-2pt}\\
\vspace{6pt}  
 \vdots   & &\vdots&\vdots&\vdots&\vdots&\vdots&   &   \\
 \uparrow & &   0  &   1  &   2  &   3  &   4  &\leftarrow\,j\\
m&&&&&&&&      
\end{array}
\]
\end{exa}

We define a partial order $\triangleleft$ on $\Z^N_+$ by 
\[
\bk'\triangleleft\bk\iff \ubj'=\ubj\ \text{and}\ \ubm'\prec\ubm.
\]

\subsubsection{}\label{eigenvec}

\begin{dfn}
For $\bk\in\Z^N_+$ 
we define $\Psi^\bk=\Phi^{{\ubm}}\otimes \bv_{{\ubj}}\in \wedge^NV(z)$.
\end{dfn}
\begin{rem}
Takemura-Uglov introduced vectors $\Phi^{{\bm}}\otimes \bv_{{\bj}}\in \wedge^NV(z)$ in \cite{takemura-uglov-level0}, which are different from ours. 
\end{rem}
\begin{prop}
$\Psi^{\bk}$ is a simultaneous eigenvector for the actions of $K^{\pm}_{i}(z)$'s ($i\in\{1,\ldots,l-1\}$). 
\end{prop}
\begin{proof}
It follows from Theorem \ref{lemma of VV} and the definition of $\Phi^{{\ubm}}$ in \ref{def of NSMP}.
\end{proof}

\subsubsection{}\label{prop2}
For $\bk\in\Z^N_+$ we define
\[
\varepsilon(\bk)=\sharp\{(a,b)\mid a<b, {\sm}(a)>{\sm}(b)\}.
\]
\begin{prop}
\[
\Psi^\bk=(-q)^{\varepsilon(\bk)}u_\bk+\sum_{\bk'\triangleleft\,\bk}c(\bk,\bk')u_{\bk'}\quad (\exists c(\bk,\bk')\in\mathcal{R}).
\]
\end{prop}
\begin{proof}
By the definition of nonsymmetric Macdonald polynomials in \ref{def of NSMP}, 
\[
\Phi^{\ubm}\otimes \bv_{\ubj}=z^{\ubm}\otimes \bv_{\ubj}+\sum_{\bm'\prec\ubm}c(\ubm,\bm')z^{\bm'}\otimes\bv_{\ubj} \quad (\exists c(\ubm,\bm')\in\mathcal{R}).
\] 
On the other hand by the relation in \ref{q-wedge} we can verify 
\begin{align*}
z^{\bm'}\otimes \bv_{\ubj}=&\,(-q)^{c(\bm',\,\ubj)}z^{\widehat{\sigma_{\bm'}}(\bm')}\otimes \bv_{\widehat{\sigma_{\bm'}}(\ubj)}\\
&+\sum_{\bm''\prec\bm'}c'(\bm',\bm'')z^{\widehat{\sigma_\bm''}(\bm'')}\otimes \bv_{\widehat{\sigma_{\bm''}}(\ubj)} \quad (\exists c'(\bm',\bm'')\in\mathcal{R})
\end{align*}
and $c(\ubm,\,\ubj)=\varepsilon(\bk)$. Then the statement follows.
\end{proof}

We define
\[
\Z^N_{+,\gamma }=\{\bk\in\Z^N_{+}\mid \ubm\in\mathcal{M}^{c,\gamma}_{N,l}\}.
\]
Note that if $\bk\in\Z^N_{+,\gamma }$ and $\bk'\triangleleft\bk$, then $\bk'\in\Z^N_{+,\gamma }$.

\begin{cor}
If $N>\gamma l$, then
$\{\Psi^\bk\mid \bk\in\Z^N_{+,\gamma }\}$ is a basis of $V^{c,\gamma}_N$.
\end{cor}

\subsubsection{}\label{eigenval}

For $\bk\in\Z^{\alpha l+c}_{+,\gamma }$ ($\alpha l+c>\gamma l$),
let $\lambda$ denote the Young diagram corresponding to $\iota^c_{\alpha l+c,\infty}(\bk)\in\Pi_c$.

\begin{prop}For $i=1,\ldots,l-1$ we have
\begin{align*}
&K_i^{\pm}(z)(\Psi^\bk)=\\
&\left(\prod_{A\in A_{\lambda,i}}\frac{(st)^{1/2}\node{A}^*t^{-c-1}z-(st)^{-1/2}}{\node{A}^*t^{-c-1}z-1}
\prod_{R\in R_{\lambda,i}}\frac{(st)^{-1/2}\node{R}^*t^{-c-1}z-(st)^{1/2}}{\node{R}^*t^{-c-1}z-1}\right)^\pm\Psi^\bk.
\end{align*}
\end{prop}

\begin{proof}
By {Theorem} \ref{lemma of VV} and the defining relations $K^\pm_i(z)=\widetilde{K^\pm_i(r^{\,l-i}z)}$, it is sufficient to show
\begin{align*}
&\Phi^{\ubm}
\prod_{j_a=i-1}\theta_1^{\pm}\left(q^{l-i+1}r^{\,l-i}Y_{a}z\right)\prod_{j_b=i}\theta_{-1}^{\pm}\left(q^{l-i-1}r^{\,l-i}Y_{b}z\right)\\
=&\,\left(\prod_{A\in A_{\lambda,i}}\frac{(st)^{1/2}\node{A}^*t^{-c}z-(st)^{-1/2}}{\node{A}^*t^{-c}z-1}
\prod_{R\in R_{\lambda,i}}\frac{(st)^{-1/2}\node{R}^*t^{-c}z-(st)^{1/2}}{\node{R}^*t^{-c}z-1}\right)^\pm\Phi^{\ubm}.
\end{align*}
First we have
\begin{align*}
\Phi^{\ubm}q^{l-i+1}r^{\,l-i}Y_{a}&=q^{l-i+1}r^{\,l-i}q^{1-N}p^{m_a}q^{\sigma_{\ubm}(a)-N-1}\Phi^{\ubm}\\
&=s^{\sigma_{\ubm}(a)-N+1/2}t^{\,\sigma_{\ubm}(a)-N-i+l(m_a+1)+1/2}\Phi^{\ubm}\\
&=s^{-\widehat{\sigma_{\bj}}^{-1}(a)+3/2}t^{-\widehat{\sigma_{\bj}}^{-1}(a)-k_a+1/2}\Phi^{\ubm},\\
\Phi^{\ubm}q^{l-i-1}r^{\,l-i}Y_{b}&=s^{\sigma_{\ubm}(b)-N-1/2}t^{\,\sigma_{\ubm}(b)-N-i+l(m_b+1)-1/2}\Phi^{\ubm}\\
&=s^{-\widehat{\sigma_{\bj}}^{-1}(b)+1/2}t^{-\widehat{\sigma_{\bj}}^{-1}(b)-k_b+1/2}\Phi^{\ubm}.
\end{align*}
We classify the elements of 
$\{a\mid j_a=i-1\}\cup\{b\mid j_b=i\}$
into three types :
\begin{enumerate}
\item[(1)] $a$ and $b$ such that $m_a=m_b$, $j_a=i-1$, $j_b=i$,
\item[(2)] $a$ such that $j_a=i-1$ and $(m_a,i)\notin\bk$, and
\item[(3)] $b$ such that $j_b=i$ and $(m_b,i-1)\notin\bk$.
\end{enumerate}
In the case of type (1), 
we have $\widehat{\sigma_{\bj}}^{-1}(a)-1=\widehat{\sigma_{\bj}}^{-1}(b)$, 
$k_a+1=k_b$. Thus 
\[
\Phi^{\ubm}q^{l-i}r^{\,l-i-1}Y_{a}=\Phi^{\ubm}q^{l-i-2}r^{\,l-i-1}Y_{b},
\]
and so 
\[
\Phi^{\ubm}\theta_1\left(q^{l-i}Y_ar^{l-i-1}z\right)\theta_{-1}\left(q^{l-i-2}Y_br^{l-i-1}z\right)=\Phi^{\ubm}.
\]
In the case of type (2), the node 
${A}=(\widehat{\sigma_{\bj}}^{-1}(a)-1,\widehat{\sigma_{\bj}}^{-1}(a)+k_a-c)$ is an addable $i$-node.
We have
\begin{align*}
\Phi^{\ubm}\theta_1\left(q^{l-i}Y_ar^{l-i-1}z\right)
&=\theta_1\left(s^{-\widehat{\sigma_{\bj}}^{-1}(a)+3/2}t^{-\widehat{\sigma_{\bj}}^{-1}(a)-k_a+1/2}z\right)\Phi^{\ubm}\\
&=\theta_1\left(s^{1/2}t^{-c+1/2}\node{A}^*z\right)\Phi^{\ubm}\\
&=\frac{(st)^{1/2}\node{A}^*t^{-c}z-(st)^{-1/2}}{\node{A}^*t^{-c}z-1}\Phi^{\ubm}.
\end{align*}
In the case of type (3), the node 
${R}=(\widehat{\sigma_{\bj}}^{-1}(b)-1,\widehat{\sigma_{\bj}}^{-1}(b)+k_b-c-1)$ is a removable $i$-node.
We have
\begin{align*}
\Phi^{\ubm}\theta_1\left(q^{l-i-2}Y_br^{l-i-1}z\right)
&=\theta_1\left(s^{-\widehat{\sigma_{\bj}}^{-1}(b)+1/2}t^{-\widehat{\sigma_{\bj}}^{-1}(b)-k_b+1/2}z\right)\Phi^{\ubm}\\
&=\theta_1\left(s^{-1/2}t^{-c-1/2}\node{R}^*z\right)\Phi^{\ubm}\\
&=\frac{(st)^{-1/2}\node{R}^*t^{-c}z-(st)^{1/2}}{\node{R}^*t^{-c}z-1}\Phi^{\ubm}.
\end{align*}
Thus the claim follows.
\end{proof}

\subsubsection{}\label{mult-free}

\begin{prop}
If $\bk,\bk'\in\Z^{\alpha l+c}_{+,\gamma }$ ($\alpha l+c>\gamma l$) and the eigenvalues of $K^\pm_i(z)$ for $\Psi^\bk$ and $\Psi^{\bk'}$ coincide for all $i\in\{1,\ldots,l-1\}$, then $\bk=\bk'$.
\end{prop}
\begin{proof}
The coincidence of the eigenvalues of $K^\pm_i(z)$ implies
\begin{align*}
&\prod_{A\in A_{\lambda,i}}(st)^{1/2}\node{A}^*t^{-c}z-(st)^{-1/2}
\prod_{R\in R_{\lambda,i}}(st)^{-1/2}\node{R}^*t^{-c}z-(st)^{1/2}\\
&\times
\prod_{A\in A_{\lambda',i}}\node{A}^*t^{-c}z-1
\prod_{R\in R_{\lambda',i}}\node{R}^*t^{-c}z-1\\
=&
\prod_{A\in A_{\lambda,i}}\node{A}^*t^{-c}z-1
\prod_{R\in R_{\lambda,i}}\node{R}^*t^{-c}z-1\\
&\times
\prod_{A\in A_{\lambda',i}}(st)^{1/2}\node{A}^*t^{-c}z-(st)^{-1/2}
\prod_{R\in R_{\lambda',i}}(st)^{-1/2}\node{R}^*t^{-c}z-(st)^{1/2}.
\end{align*}
Since $\left|\{(s,t)\mid s-t=n\}\cap\left(A_{\lambda}\cup R_\lambda\right)\right|<1$ for any $n\in\Z$, we have
\begin{align*}
&\prod_{A\in A_{\lambda,i}}\left((st)^{1/2}\node{A}^*t^{-c}z-(st)^{-1/2}\right)\\
&\times\prod_{R\in R_{\lambda,i}}\left((st)^{-1/2}\node{R}^*t^{-c}z-(st)^{1/2}\right)\bigg|_{z=\node{X}t^c}\neq 0
\end{align*} 
for any $X\in A_{\lambda,i}\cup R_{\lambda,i}$.
So we have 
$X\in A_{\lambda',i}\cup R_{\lambda',i}$,
and it follows that $A_{\lambda,i}\cup R_{\lambda,i}=A_{\lambda',i}\cup R_{\lambda',i}$.

It is easy to see the set $\bigcup_{i\neq 0}\left(A_{\lambda,i}\cup R_{\lambda,i}\right)$ determines $\lambda$. So the claim follows.
\end{proof}

\begin{cor}
If $X\in\wedge^NV(z)$ is a simultaneous eigenvector for the actions of $K^{\pm}_{i}(z)$'s ($i\in\{1,\ldots,l-1\}$) and 
\[
X=(-q)^{\varepsilon(\bk)}u_\bk+\sum_{\bk'\triangleleft\,\bk}c(\bk')u_{\bk'}\quad (\exists c(\bk')\in\mathcal{R}),
\]
for $\bk\in\Z^{\alpha l+c}_{+,\gamma }$ ($\alpha l+c>\gamma l$),
then $X=\Psi^\bk$.
\end{cor}
\begin{proof}
It follows from Corollary \ref{prop2} and the previous proposition.
\end{proof}

\subsubsection{}\label{prop3}


\begin{prop}
For $\bk\in\Z^{\alpha l+c}_{+,\gamma }$ ($\alpha l+c>\gamma l$), $\Psi^\bk$ is a simultaneous eigenvector for the actions of $K^\pm_i(z)$'s ($i\in\{0,\ldots,l-1\}$).
\end{prop}
\begin{proof}
By definition $K^\pm_{i}(z)$'s commute with each other.
Notice that a matrix which commutes with a diagonal matrix with diagonal elements different from each other is diagonal.
It follows from Proposition \ref{eigenvec}, Corollary \ref{prop2} and Proposition \ref{mult-free} 
that $\Psi^\bk$ is also a simultaneous eigenvector for the action of $K^\pm_0(z)$.
\end{proof}

\subsubsection{}\label{hoge}

\begin{prop}
For $\bk\in\Z^{al+c}_{+,\gamma }$ ($al+c>\gamma l$), we have
\[
\rho(q^{-\varepsilon(\bk)}\Psi^{\bk})=q^{-\varepsilon(\rho(\bk))}\Psi^{\rho(\bk)}.
\]
\end{prop}
\begin{proof}
By Lemma \ref{VVschur} and Theorem \ref{VVschur}, $\rho(\Psi^\bk)$ is also a simultaneous eigenvector of $K^\pm_{i}(z)$'s ($i\in\{0,\ldots,l-1\}$).
Note that 
\begin{align*}
\rho(q^{-\varepsilon(\bk)}\Psi^\bk)&=
\rho(u_\bk)+\sum_{\bk'\triangleleft\,\bk}c(\bk,\bk')\rho(u_{\bk'})\\
&=u_{\rho(\bk)}+\sum_{\bk'\triangleleft\,\bk}c(\bk,\bk')u_{\rho(\bk')}
\end{align*}
and $\rho$ preserves the order $\triangleleft$.
Then the statement follows from Corollary \ref{mult-free}.
\end{proof}

\begin{cor}
For $\bk\in\Z^{\alpha l+c}_{+,\gamma }$ ($\alpha l+c>\gamma l$), the eigenvalue of $\Psi^\bk$ for $K^\pm_{0}(z)$ is given by the same formula as in Proposition \ref{eigenval}. 
\end{cor}

\subsubsection{}\label{lift}

For $\beta>\alpha$ we write simply $\iota^{\,c}_{\alpha,\beta}$ for $\iota^{\,c}_{\alpha l+c,\beta l+c}$. 

\begin{lem}
For $\bk\in\Z^{\alpha l+c}_{+,\gamma }$ ($\alpha l+c>\gamma l$), we have
\[
\iota^{\,c}_{\alpha,\beta}(q^{-\varepsilon(\bk)}\Psi^\bk)=q^{-\varepsilon(\iota^{\,c}_{\alpha,\beta}(\bk))}\Psi^{\iota^{\,c}_{\alpha,\beta}(\bk)}.
\]
\end{lem}
\begin{proof}
By {Theorem} \ref{TU}, $\iota^{\,c}_{\alpha,\beta}(\Psi^\bk)\in\wedge^{\beta l+c}V(z)$ is also a simultaneous eigenvector.
Note that
\begin{align*}
\iota^{\,c}_{\alpha,\beta}(q^{-\varepsilon(\bk)}\Psi^\bk)&=
\iota^{\,c}_{\alpha,\beta}(u_{\bk})+\sum_{\bk'\triangleleft\,\bk}c(\bk,\bk')\iota^{\,c}_{\alpha,\beta}(u_{\bk'})\quad (\exists c(\bk,\bk')\in\mathcal{R})\\
&=u_{\iota^{\,c}_{\alpha,\beta}(\bk)}+\sum_{\bk'\triangleleft\,\bk}c(\bk,\bk')u_{\iota^{\,c}_{\alpha,\beta}(\bk')},
\end{align*}
and $\iota^{\,c}_{\alpha,\beta}$ preserves the order $\triangleleft$.
The claim follows from Corollary \ref{mult-free}.
\end{proof}
This lemma assures us of the well-definedness for the following definition :
\begin{dfn}
For a Young diagram $\lambda$, take $\bk\in\Z^{\alpha l+c}_{+,\gamma }$ ($\alpha l+c>\gamma l$) such that
$\iota^{\,c}_{\alpha l+c,\infty}(\bk)\in\Pi_c$ corresponds to $\lambda$. We define
\[
\Psi^\lambda_c=\iota^{\,c}_{\alpha l+c,\infty}(q^{-\varepsilon(\bk)}\Psi^\bk)\in F(c).
\]
\end{dfn}

\section{Isomorphism}\label{isom}

{\ }

\vspace{-14pt}

In this section we construct an isomorphism between the equivariant K-groups of the quiver varieties and the q-Fock space.

First we express the actions of $e_{i.n}$'s on $\Psi^\bk$, which can be done completely in terms of Young diagrams (Theorem \ref{formula}).
This is basically due to Proposition \ref{prop of NSMP} although we need the residue theorem and a little complicated induction.

After suitable renormalizations, we arrive at the isomorphism.

\subsection{Formula for the representation on the q-Fock space}

\subsubsection{}\label{key lemma}

We will give a formula for the action of $e_{i,n}$ on $\Psi^\bk$.
For this the following proposition is essential ;
\begin{prop}
\begin{enumerate}
\item[(1)]
For $m_1<\cdots<m_N$ and $j_1<j_2=\cdots=j_N$ 
we have
\[
\Phi^\bm\left(\sum_{a=1}^{N}T_{a,1}\right)\otimes \bv_\bj
=\sum_{a=1}^N
\left(\prod_{b=a+1}^{N}\frac{-\zeta_{b}(\bm)+q^2\zeta_{a}(\bm)}{\zeta_{b}(\bm)-\zeta_{a}(\bm)}\right)\Phi^{\bm(a)}\otimes \bv_{\bj}
\]
where $\bm(a)=(m_a,m_1,\ldots,\overset{\wedge}{m_a},\ldots)$.
\item[(2)]
For $m_1<\cdots<m_{N-1}$, $m_N=m_a$ ($1\leq a<N$) and $j_1=j_2=\cdots=j_N$ we have
\[
\Phi^\bm\otimes \bv_\bj=0
\]
\item[(3)]
For $\bk$ such that 
$m_1<\cdots<m_{N-1}$, $m_a<m_N<m_{a+1}$ and $j_1=j_2=\cdots=j_N$ we have
\[
\Phi^\bm\otimes \bv_\bj
=
\left(\prod_{b=a+1}^{N-1}\frac{-\zeta_{b}(\bn)+q^2\zeta_{a}(\bn)}{\zeta_{b}(\bn)-\zeta_{a}(\bn)}\right)\Phi^{\overline{\bm}}\otimes \bv_{\bj}
\]
where $\overline{\bm}=(\ldots,m_a,m_N,m_{a+1},\ldots)$.
\end{enumerate}
\end{prop}

\begin{proof}
If $j_a=j_{a+1}$ and $m_a>m_{a+1}$ we have
\begin{align*}
q^2\Phi^\bm\otimes\bv_\bj&=\Phi^\bm\otimes\vt_a\bv_\bj\quad &\text{by (\ref{left action})}\\
&=\Phi^\bm\pt_a\otimes\bv_\bj\quad &\\
&=\frac{(-q^2+1)}{x-1}\Phi^\bm\otimes\bv_\bj-\frac{(x-q^2)(q^2x-1)}{(x-1)^2}\Phi^{\sigma_a\bm}\otimes\bv_\bj &\text{by (\ref{prop of NSMP})}.
\end{align*}
where $\sigma_a\bm=(\ldots,m_{a+1},m_a,\ldots)$ and $x=\frac{\zeta_{a+1}(\bm)}{\zeta_a(\bm)}$. Thus
\[
\Phi^\bm\otimes\bv_\bj=\frac{-\zeta_{a+1}(\bm)+q^2\zeta_{a}(\bm)}{\zeta_{a+1}(\bm)-\zeta_{a}(\bm)}\Phi^{\sigma_a\bm}\otimes\bv_\bj.
\]
The statement of $(3)$ follows this.

If $j_a=j_{a+1}$ and $m_a=m_{a+1}$ then $x=\frac{\zeta_{a+1}(\bm)}{\zeta_a(\bm)}=q^2$.
Thus we have
\begin{align*}
q^2\Phi^\bm\otimes\bv_\bj&=\Phi^\bm\otimes\vt_a\bv_\bj\quad &\text{by (\ref{left action})}\\
&=\Phi^\bm\pt_a\otimes\bv_\bj\quad &\\
&=\frac{(-q^2+1)}{x-1}\Phi^\bm\otimes\bv_\bj &\text{by (\ref{prop of NSMP})}\\
&=-\Phi^\bm\otimes\bv_\bj,
\end{align*}
and so $\Phi^\bm\otimes\bv_\bj=0$. This shows $(2)$.

We will prove $(1)$ by induction for $N$.
Assume the statement is true for $N'<N$.
Then
\begin{align*}
\Phi^\bm\left(\sum_{a=1}^{N}T_{a,1}\right)\otimes \bv_\bj&=
\Phi^\bm\left(\left(\sum_{a=2}^{N}T_{a,2}\right)\vt_1+1\right)\otimes \bv_\bj\\
&=\left[\sum_{a=2}^{N}
\left(\prod_{b=a+1}^{N}\frac{-\zeta_{b}(\bk)+q^2\zeta_{a}(\bk)}{\zeta_{b}(\bk)-\zeta_{a}(\bk)}\right)\Phi^{\widetilde{\bm(a)}}\vt_1+\Phi^{\bm}\right]\otimes \bv_{\bj},
\end{align*}
where $\widetilde{\bm(a)}=\bm(a)=(m_1,m_a,m_2,\ldots,\overset{\wedge}{m_a},\ldots)$.
Here we use the assumption of induction.
Although the situations are not exactly same, 
commutativity of $\vt_1$ with $\vt_a\ (a\geq 3)$ allows us a parallel argument. 

Further we have
\begin{align*}
\Phi^{\widetilde{\bm(a)}}\vt_1\otimes \bv_{\bj}&=\frac{(-q^2+1)}{x-1}\Phi^{\widetilde{\bm(a)}}\otimes \bv_{\bj}-\Phi^{\sigma_1\widetilde{\bm(a)}}\otimes \bv_{\bj}\\
&=\frac{(-q^2+1)\zeta_1(\bm)}{\zeta_{a}(\bm)-\zeta_1(\bm)}\prod_{b=2}^{a-1}\frac{-\zeta_{b}(\bm)+q^2\zeta_{a}(\bm)}{\zeta_{b}(\bm)-\zeta_{a}(\bm)}\Phi^{\bm}\otimes \bv_{\bj}-\Phi^{\bm(a)}\otimes \bv_{\bj}.
\end{align*}
We can see the coefficients of $\bm(a)\ (a\geq 2)$ coincide with required ones.
For the coefficient of $\bm$ we need to check

\[
\left(\prod_{b=2}^{N}\frac{\zeta_{b}(\bk)-q^2\zeta_{1}(\bk)}{\zeta_{b}(\bk)-\zeta_{1}(\bk)}\right)
=\sum_{a=2}^{N}
\left(
\frac{(-q^2+1)\zeta_1(\bm)}{\zeta_{a}(\bm)-\zeta_1(\bm)}
\prod_{b\neq a}\frac{\zeta_{b}(\bk)-q^2\zeta_{a}(\bk)}{\zeta_{b}(\bk)-\zeta_{a}(\bk)}
\right)+1.
\]
This follows the next lemma.
\end{proof}

\begin{lem}
\[
\sum_{a=1}^N\left(\frac{-q^2+1}{x_1-1}\prod_{b\neq a}\frac{x_b-q^2x_a}{x_b-x_a}\right)
=\left(\prod_{a=1}^N\frac{x_a-q^2}{x_a-1}\right)-1
\]
\end{lem}
\begin{proof}
Apply the residue theorem for a rational function 
\[
f(Z)=\frac{1}{Z(Z-1)}\prod_{a=1}^N\frac{x_a-q^2Z}{x_a-Z}.
\]
\end{proof}

\subsubsection{}\label{formula}

\begin{thm}
\begin{align*}
e_{i,n}(\Psi_c^\lambda)&=
\sum_{X\in R_{\lambda,i}}\hspace{2pt}\left(t^c\node{X}\right)^{n}\hspace{-2pt}
\,\left(\,
\prod_{
\begin{subarray}{c}
A\in A_\lambda,i\\
A>X
\end{subarray}
}
\frac{-(st)^{-1/2}\node{X}^{*}+(st)^{1/2}\node{A}^{*}}{\node{X}^{*}-\node{A}^{*}}
\right)\\
&\hspace{22mm}\times\left(\,
\prod_{
\begin{subarray}{c}
R\in R_{\lambda\backslash X,i}\\
R>X
\end{subarray}
}
\frac{-(st)^{-1/2}\node{R}^{*}+(st)^{1/2}\node{X}^{*}}{\node{R}^{*}-\node{X}^{*}}
\right)
\,\Psi_c^{\lambda\backslash X},\\
f_{i,n}(\Psi_c^\lambda)&=
\sum_{X\in A_{\lambda,i}}\hspace{2pt}\left(t^c\node{X}\right)^{n}\hspace{-2pt}
\,\left(\,
\prod_{
\begin{subarray}{c}
A\in A_{\lambda\cup X},i\\
A<X
\end{subarray}
}
\frac{-(st)^{-1/2}\node{X}^{*}+(st)^{1/2}\node{A}^{*}}{\node{X}^{*}-\node{A}^{*}}
\right)\\
&\hspace{22mm}\times\left(\,
\prod_{
\begin{subarray}{c}
R\in R_{\lambda,i}\\
R<X
\end{subarray}
}
\frac{-(st)^{-1/2}\node{R}^{*}+(st)^{1/2}\node{X}^{*}}{\node{R}^{*}-\node{X}^{*}}
\right)
\,\Psi_c^{\lambda\cup X},\\
K^\pm_i(z)(\Psi_c^\lambda)&=
\left(\prod_{A\in A_{\lambda,i}}\frac{(st)^{1/2}\node{A}^*t^{-c}z-(st)^{-1/2}}{\node{A}^*t^{-c}z-1}
\prod_{R\in R_{\lambda,i}}\frac{(st)^{-1/2}\node{R}^*t^{-c}z-(st)^{1/2}}{\node{R}^*t^{-c}z-1}\right)^\pm\Psi_c^\lambda.
\end{align*}
\end{thm}

\begin{proof}
The formulas for $K^\pm_i(z)$'s are nothing but Proposition \ref{eigenval} and Corollary \ref{hoge}.
We will check for $e_{i,n}$'s. 
 
For $a,b\in \{1,\ldots,N\}$
we put 
\[
f(a,b)=\frac{q^{-1}\zeta_b(\ubm)-q\zeta_a(\ubm)}{\zeta_b(\ubm)-\zeta_a(\ubm)}.
\]
Then we have

\begin{eqnarray*}
&e_{i,n}&\hspace{-5.3mm}\left(\Phi^\ubm\otimes\bv_\ubj\right)\\
&\underset{\ref{lemma of VV}}{=}&q^{1-n_{i}}\Phi^\ubm\left(\sum_{a=\bar{n}_{i-1}+1}^{\bar{n}_i}T_{a,\bar{n}_{i-1}+1}\right)
\left(q^{l-i}Y_{\bar{n}_{i-1}+1}r^{l-i}\right)^{-n}
\otimes \bv_{\ubj^-}\\
&\underset{\ref{key lemma} (1)}{=}&q^{1-n_{i}}\sum_{a=\bar{n}_{i-1}+1}^{\bar{n}_i}\left(\prod_{b=a+1}^{\bar{n}_i}-qf(a,b)\right)\Phi^{\ubm(a)}
\left(q^{l-i}Y_{\bar{n}_{i-1}+1}r^{l-i}\right)^{-n}
\otimes \bv_{\ubj^-}\\
&\underset{\text{see \ref{eigenval}}}{=}&q^{1-n_{i}}\sum_{a=\bar{n}_{i-1}+1}^{\bar{n}_i}\Bigl(t^c\node{X_a}\Bigr)^{n}\left(\prod_{b=a+1}^{\bar{n}_i}-qf(a,b)\right)\Phi^{\ubm(a)}
\otimes \bv_{\ubj^-}\\
&\underset{\ref{key lemma} (2)}{=}&q^{1-n_{i}}\sum_{
\begin{subarray}{c}
\bar{n}_{i-1}<a\leq\bar{n}_i\\
(\underline{m}_a,i-1)\notin\bk
\end{subarray}
}
\Bigl(t^c\node{X_a}\Bigr)^{n}
\left(\prod_{b=a+1}^{\bar{n}_i}-qf(a,b)\right)\left(\prod_{b=a'+1}^{\bar{n}_{i-1}}-qf(b,a)\right)\Phi^{\overline{\ubm(a)}}\otimes \bv_{\ubj^-}\\
&=&\hspace{-6mm}\sum_{
\begin{subarray}{c}
\bar{n}_{i-1}<a\leq\bar{n}_i\\
(\um{a},i-1)\notin\bk
\end{subarray}
}
q^{n_{i-1}^+(a)-n_{i}^-(a)}\Bigl(t^c\node{X_a}\Bigr)^{n}\left(\prod_{b=a+1}^{\bar{n}_i}-f(a,b)\right)\left(\prod_{b=a'+1}^{\bar{n}_{i-1}}-f(b,a)\right)\Phi^{\overline{\ubm(a)}}\otimes \bv_{\ubj^-}, 
\end{eqnarray*}
where 
\begin{itemize}
\item
$\ubm(a)=(\ldots,\um{\bar{n}_i},\um{a},\um{\bar{n}_i+1},\ldots,\overset{\wedge}{\um{a}},\ldots)$, \\
$\overline{\ubm(a)}=(\ldots,\um{a'},\um{a},\um{a'+1},\ldots,\overset{\wedge}{\um{a}},\ldots)$,
\item $X_a=\left(\sj^{-1}(a)-1,\sj^{-1}(a)+k_a-c\right)$ denote the top node on the $a$-th line of $\lambda$, and 
\item
$n_{i-1}^+(a)=\bar{n}_i-a'$, $n_{i}^-(a)=a-\bar{n}_i-1$.
\end{itemize}
Since $\varepsilon(\ubm,\,\ubj)-\varepsilon(\overline{\ubm(a)},\,\ubj^-)=n_{i-1}^+(a)-n_{i}^-(a)$ we have
\begin{align*}
&e_{i,n}\left(q^{-\varepsilon(\ubm,\,\ubj)}\Phi^{\ubm}\otimes\bv_{\ubj}\right)
=\\
&\sum_{
\begin{subarray}{c}
\bar{n}_{i-1}<a\leq\bar{n}_i\\
(\um{a},i-1)\notin\bk
\end{subarray}
}
\Bigl(t^c\node{X_a}\Bigr)^{n}
\left(\prod_{b=a+1}^{\bar{n}_i}-f(a,b)\right)\left(\prod_{b=a'+1}^{\bar{n}_{i-1}}-f(b,a)\right)
q^{-\varepsilon(\overline{\ubm(a)},\,\ubj^-)}\Phi^{\overline{\ubm(a)}}\otimes \bv_{\ubj^-}.
\end{align*}
As in the proof of Proposition \ref{eigenval} we can arrange the right hand side of the above equation 
by classify the element of $\{a+1,\ldots,\bar{n}_i\}\cup\{a'+1,\ldots,\bar{n}_{i-1}\}$ into three types, 
and finally we get  
\begin{align*}
e_{i,n}(\Psi_c^\lambda)=
\sum_{X\in R_{\lambda,i}}\hspace{2pt}\Bigl(t^c\node{X}\Bigr)^{n}\hspace{-2pt}
&\,\left(\,
\prod_{
\begin{subarray}{c}
A\in A_\lambda,i\\
A>X
\end{subarray}
}
\frac{-(st)^{-1/2}\node{X}^{*}+(st)^{1/2}\node{A}^{*}}{\node{X}^{*}-\node{A}^{*}}
\right)\\
&\ \ \left(\,
\prod_{
\begin{subarray}{c}
R\in R_{\lambda\backslash X,i}\\
R>X
\end{subarray}
}
\frac{-(st)^{-1/2}\node{R}^{*}+(st)^{1/2}\node{X}^{*}}{\node{R}^{*}-\node{X}^{*}}
\right)
\,\Psi_c^{\lambda\backslash X}.
\end{align*}

\end{proof}

\subsection{Normalizations}

\subsubsection{}\label{normalization1}

\begin{dfn}
For $\lambda\in\Pi$ we define 
\[
N(\lambda)=\prod\left(s^{x_h-x_t}t^{y_h-y_t}-1\right),
\]
where the product runs over all $l$-hooks $((x_h,y_h),(x_t,y_t))$. 
\end{dfn}
We can easily verify the following lemma : 
\begin{lem}
If $X$ is a removable $i$-node of $\lambda$, then we have 
\begin{align*}
N(\lambda)/N(\lambda\backslash X)=&
\prod_{
\begin{subarray}{c}
A\in A_{\lambda,i}\\
A<X
\end{subarray}
}
\left(\node{A}\node{X}^{*}-1\right)
\prod_{
\begin{subarray}{c}
A\in A_{\lambda,i}\\
A>X
\end{subarray}
}
\left(st\,\node{X}\node{A}^{*}-1\right)\\
&\prod_{
\begin{subarray}{c}
R\in R_{\lambda,i}\\
R<X
\end{subarray}
}
\left(st\,\node{R}\node{X}^{*}-1\right)^{-1}
\prod_{
\begin{subarray}{c}
R\in R_{\lambda,i}\\
R>X
\end{subarray}
}
\left(\node{X}\node{R}^{*}-1\right)^{-1}.
\end{align*}
\end{lem}

\begin{rem}
From geometrical point of view, $N(\lambda)$ is derived from the Kozsul complex of the unstable manifold 
, with respect to a specific $\C^*$-action, on which points converge to the fixed point $\lambda$.
\end{rem}

\subsubsection{}\label{normalization2}

For $\mu\in\Pi$ we will define $M(\mu)\in\R$ inductively.
First we set $M(\emptyset)=1$. 
Let $Y=(a,b)$ be the most right node of the top row of $\mu$.
Then we set
\begin{align*}
M(\mu)=&\,
M(\mu\backslash Y)
\left(s^{-1/2}t^{1/2}\right)^{\alp_{i-1}(\mu)}
\node{Y}^{\,\delta(b\equiv 0)}\\
&\times\prod_{
\begin{subarray}{c}
A\in A_{\mu,i}\\
A<Y
\end{subarray}
}
\left((st)^{-1/2}\node{A}\right)
\prod_{
\begin{subarray}{c}
R\in R_{\mu\backslash Y,i}\\
R<Y
\end{subarray}
}
\left((st)^{-1/2}\node{R}^{*}\right).
\end{align*}


\begin{lem}
If $X$ is a removable $i$-node of $\lambda\in\Pi$, then we have
\begin{align*}
M(\lambda)=&\,M(\lambda\backslash X)
\left(s^{-1/2}t^{1/2}\right)^{\alp_{i-1}(\lambda)}\\
&\times\prod_{
\begin{subarray}{c}
A\in A_{\lambda,i}\\
A<X
\end{subarray}
}
\left((st)^{-1/2}\node{A}\right)
\ \prod_{
\begin{subarray}{c}
A\in A_{\lambda,i}\\
A>X
\end{subarray}
}
\node{X}\\
&\times\prod_{
\begin{subarray}{c}
R\in R_{\lambda\backslash X,i}\\
R<X
\end{subarray}
}
\left((st)^{-1/2}\node{R}^{*}\right)
\prod_{
\begin{subarray}{c}
R\in R_{\lambda\backslash X,i}\\
R>X
\end{subarray}
}
\node{X}^{*}.
\end{align*}
\end{lem}

\begin{proof}
We divide $\lambda$ into $\lambda_l$ and $\lambda_r$ by the vertical line on the right of $X$.

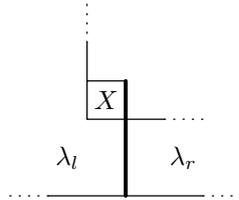
\begin{center}
\input{epsilon2.tpc}
\end{center}

Let us write $\mu<\lambda$ if we can get $\mu$ from $\lambda$ by successive removing the nodes on the top of the most right line of diagrams. 
Take $\mu\in\Pi$ such that $\lambda_l\lneq\mu<\lambda$.
Let $Y$ be the node on the top of the most right line of $\mu$.

Then we can verify
\begin{align*}
&\frac{M(\mu)}{M(\mu\backslash Y)}\bigg/\frac{M(\mu\backslash X)}{M(\mu\backslash (X\cup Y))}\\
=&\,\left(s^{-1/2}t^{1/2}\right)^{\delta(j-1\equiv i)}
\left((st)^{-1/2}\node{X}^{*}\right)^{\delta(j\equiv i)}
\left((st)^{-1/2}\node{X}\right)^{-\delta(j\equiv i)}\\
&\times\left((st)^{-1/2}\node{X}t\right)^{\delta(j\equiv i-1)}
\left((st)^{-1/2}\node{X}s\right)^{\delta(j\equiv i+1)}\\
=&\,
\left(s^{-1/2}t^{1/2}\right)^{\delta(j-1\equiv i)}
\node{X}^{\,\delta(j\equiv i-1)-2\delta(j\equiv i)+\delta(j\equiv i+1)}.
\end{align*}
So we have
\begin{align*}
&\frac{M(\lambda)}{M(\lambda_l)}\bigg/\frac{M(\lambda\backslash X)}{M(\lambda_l\backslash X)}\\
=&\,\left(s^{-1/2}t^{1/2}\right)^{\alp_{i-1}(\lambda_r)}
\node{X}^{\,\alpha_{i-1}(\lambda_r)-2\alpha_{i}(\lambda_r)+\alpha_{i+1}(\lambda_r)}.
\end{align*}
Let $j$ denote the content of the node on the bottom of the most left line be $j$ of $\lambda_r$. 
Note that we have $\delta(b\equiv 0)=\delta(j\equiv i)$ and 
\[
|A_{\lambda_r,i}|-|R_{\lambda_r,i}|=\alpha_{i-1}(\lambda_r)-2\alpha_{i}(\lambda_r)+\alpha_{i+1}(\lambda_r)+\delta(j\equiv i).
\]
Finally we have
\begin{align*}
&\frac{M(\lambda)}{M(\lambda\backslash X)}\\
=&\,\left(s^{-1/2}t^{1/2}\right)^{\alp_{i-1}(\lambda_r)}
\node{X}^{\,\alpha_{i-1}(\lambda_r)-2\alpha_{i}(\lambda_r)+\alpha_{i+1}(\lambda_r)}\\
&\times\left(s^{-1/2}t^{1/2}\right)^{\alp_{i-1}(\lambda_l)}
\node{X}^{\,\delta(b\equiv 0)}
\prod_{
\begin{subarray}{c}
A\in A_{\lambda_l,i}\\
A<X
\end{subarray}
}
\left((st)^{-1/2}\node{A}\right)
\prod_{
\begin{subarray}{c}
R\in R_{\lambda_l\backslash X,i}\\
R<X
\end{subarray}
}
\left((st)^{-1/2}\node{R}^{*}\right)\\
=&\,
\left(s^{-1/2}t^{1/2}\right)^{\alp_{i-1}(\lambda)}
\node{X}^{\,|A_{\lambda_r,i}|-|R_{\lambda_r,i}|}
\prod_{
\begin{subarray}{c}
A\in A_{\lambda,i}\\
A<X
\end{subarray}
}
\left((st)^{-1/2}\node{A}\right)
\prod_{
\begin{subarray}{c}
R\in R_{\lambda\backslash X,i}\\
R<X
\end{subarray}
}
\left((st)^{-1/2}\node{R}^{*}\right).
\end{align*}
So the claim follows.
\end{proof}

\subsection{Main theorem}\label{main}

Now we arrive at the main theorem ;
\begin{thm}
$K^{T}_{\mathcal{R}}(\M)$ and $F(0)$ is isomorphic as representations of $\tor$.
The isomorphism is given by 
\[
N(\lambda)b_\lambda\longmapsto M(\lambda)\Psi_0^\lambda.
\] 
\end{thm}
\begin{proof}
This follows from Theorem \ref{formula of toroidal action}, Theorem \ref{formula}, Lemma \ref{normalization1} and Lmma \ref{normalization2}.
\end{proof}

\providecommand{\bysame}{\leavevmode\hbox to3em{\hrulefill}\thinspace}
\providecommand{\MR}{\relax\ifhmode\unskip\space\fi MR }
\providecommand{\MRhref}[2]{%
  \href{http://www.ams.org/mathscinet-getitem?mr=#1}{#2}
}
\providecommand{\href}[2]{#2}

{\tt
\noindent Kentaro Nagao

\noindent Department of Mathematics, Kyoto University, Kyoto 606-8502, Japan

\noindent kentaron@math.kyoto-u.ac.jp
}

\end{document}